\documentclass[11pt,letterpaper]{amsart}
\usepackage{amsfonts,amssymb,amscd,amsthm,amsmath,color,comment}
\usepackage{tikz}
\usepackage[alphabetic]{amsrefs}
\usepackage{mathrsfs,ulem}
\usepackage{times}
\makeatletter
\def\makeop#1{\expandafter\def\csname#1\endcsname{\mathop{\rm #1}\nolimits}\ignorespaces}
\def\makeoplist#1 {\def\@@tmpa{#1}\def\@@tmpb{***}%
  \ifx\@@tmpa\@@tmpb\else\makeop{#1}\expandafter\makeoplist\fi}
\makeoplist
ad Ad Alb Alt Aut Br can Card card Char Cl Coker coker Cond cond Corr
depth diag Diff Disc disc Div div dom End Ext Fil Fr Frac Frob Gal GL
Gr gr GSp GSpin Hom hom Id id Im im Ind ind Int inv Irr irr Isom
Ker ker len length Lie Nr Nrd O Ob ob ord PGL Pic Pr pr Proj PSL rank
RD Res Sgn sgn SL SO Sp Sp sp Spec Spf Spin
St st Stab stab Std std SU Supp supp Swan Sym Tor tor tors
torsion Tr tr Trd U val Vol vol ***
\def\makermlist#1 {\def\@@tmpa{#1}\def\@@tmpb{***}\ifx\@@tmpa\@@tmpb
  \else\expandafter\def\csname#1\endcsname{{\rm#1}}\expandafter\makermlist\fi}
\makermlist ab alg arith can geom opp op red sep sh ***
\makeatother
\DeclareMathAlphabet\eusm{U}{eus}{m}{n}
\def\makebb#1{\expandafter\def\csname bb#1\endcsname{{\mathbb{#1}}}\ignorespaces}
\def\makerm#1{\expandafter\def\csname rm#1\endcsname{{\rm #1}}\ignorespaces}
\def\makebf#1{\expandafter\def\csname bf#1\endcsname{{\bf #1}}\ignorespaces}
\def\makegr#1{\expandafter\def\csname gr#1\endcsname{{\mathfrak{#1}}}\ignorespaces}
\def\makescr#1{\expandafter\def\csname scr#1\endcsname{{\mathscr{#1}}}\ignorespaces}
\def\makecal#1{\expandafter\def\csname cal#1\endcsname{{\cal #1}}\ignorespaces}
\def\makeudl#1{\expandafter\def\csname udl#1\endcsname{{\underline{#1}}}\ignorespaces}
\def\doLetters#1{%
  #1A #1B #1C #1D #1E #1F #1G #1H #1I #1J #1K #1L #1M
  #1N #1O #1P #1Q #1R #1S #1T #1U #1V #1W #1X #1Y #1Z}
\def\doletters#1{%
  #1a #1b #1c #1d #1e #1f #1g #1h #1i #1j #1k #1l #1m
  #1n #1o #1p #1q #1r #1s #1t #1u #1v #1w #1x #1y #1z}
\doLetters\makebb\doLetters\makecal\doLetters\makerm\doLetters\makebf\doLetters\makescr
\doletters\makerm \doletters\makebf \doLetters\makegr \doletters\makegr \doLetters\makeudl \doletters\makeudl
    \def\setminus{\smallsetminus}
   
\def\ringO{{\scrO}}

\makeatletter   \newdimen\mina@@\mina@@=18pt
\newcommand{\xrtarw}[2][]{\mathrel{\mathop{\,\setbox\z@\vbox{\m@th
  \hbox{$\scriptstyle\;{#1}\;\;$}\hbox{$\m@th\scriptstyle\;{#2}\;\;$}}%
  \hbox to\ifdim\wd\z@>\mina@@\wd\z@\else\mina@@\fi{\rightarrowfill@
  \displaystyle}\,}\limits^{#2}\@ifnotempty{#1}{_{#1}}}}
\newcommand{\xltarw}[2][]{\mathrel{\mathop{\,\setbox\z@\vbox{\m@th
  \hbox{$\scriptstyle\;\;{#1}\;$}\hbox{$\m@th\scriptstyle\;\;{#2}\;\;$}}%
  \hbox to\ifdim\wd\z@>\mina@@\wd\z@\else\mina@@\fi{\leftarrowfill@
  \displaystyle}\,}\limits^{#2}\@ifnotempty{#1}{_{#1}}}}
\makeatother
\def\XYmatrix{\xymatrix@M=5pt} 
\def\ncmd{\newcommand}
\ncmd{\xysubset}[1][r]{\ar@<-2.5pt>@{^(-}[#1]\ar@<2.5pt>@{_(-}[#1]}
\ncmd{\XYmatrixc}[1]{\vcenter{\XYmatrix{#1}}}
\ncmd{\xyto}[1][r]{\ar@{->}[#1]}      \ncmd{\xyinj}[1][r]{\ar@{^(->}[#1]}
\ncmd{\xysurj}[1][r]{\ar@{->>}[#1]}   \ncmd{\xyline}[1][r]{\ar@{-}[#1]}
\ncmd{\xydotsto}[1][r]{\ar@{.>}[#1]}  \ncmd{\xydots}[1][r]{\ar@{.}[#1]}
\ncmd{\xyleadsto}[1][r]{\ar@{~>}[#1]} \ncmd{\xyeq}[1][r]{\ar@{=}[#1]}
\ncmd{\xyequal}[1][r]{\ar@{=}[#1]}    \ncmd{\xyequals}[1][r]{\ar@{=}[#1]}
\ncmd{\xymapsto}[1][r]{\ar@{|->}[#1]}\ncmd{\xyimplies}[1][r]{\ar@{=>}[#1]}
\ncmd{\xytofrom}[1][r]{\ar@{<->}[#1]} 
  
\def\XYTOTO[#1]^#2_#3{\xyto[#1]<0.5ex>^{#2}\xyto[#1]<-0.5ex>_{#3}}
\ncmd{\xytoto}[1][r]{\XYTOTO[#1]}
\def\beginmat{\begin{pmatrix}}\def\endmat{\end{pmatrix}}
\def\pmat#1]{{\def\beginmat{\begin{pmatrix}}\def\endmat{\end{pmatrix}}\mat#1]}}
\def\bmat#1]{{\def\beginmat{\begin{bmatrix}}\def\endmat{\end{bmatrix}}\mat#1]}}
\def\Bmat#1]{{\def\beginmat{\begin{Bmatrix}}\def\endmat{\end{Bmatrix}}\mat#1]}}
\def\vmat#1]{{\def\beginmat{\begin{vmatrix}}\def\endmat{\end{vmatrix}}\mat#1]}}
\def\Vmat#1]{{\def\beginmat{\begin{Vmatrix}}\def\endmat{\end{Vmatrix}}\mat#1]}}
\def\smat#1]{{\def\beginmat{\begin{smallmatrix}}%
  \def\endmat{\end{smallmatrix}}\left(\mat#1]\right)}}
\def\mat#1#2]{\ifcase#1\or \matA#2]\or \matAA#2]\or \matAAA#2]\fi}
\def\matA  #1#2]{\ifcase#1\or \matAB  #2]\or \matABB  #2]\or \matABBB  #2]\fi}
\def\matAA #1#2]{\ifcase#1\or \matAAB #2]\or \matAABB #2]\or \matAABBB #2]\fi}
\def\matAAA#1#2]{\ifcase#1\or \matAAAB#2]\or \matAAABB#2]\or \matAAABBB#2]\fi}
\def\matAB[#1]{\beginmat#1\endmat}
\def\matABB[#1,#2]{\beginmat#1&#2\endmat}
\def\matABBB[#1,#2,#3]{\beginmat#1&#2&#3\endmat}
\def\matAAB[#1;#2]{\beginmat#1\\#2\endmat}
\def\matAABB[#1,#2;#3,#4]{\beginmat#1&#2\\#3&#4\endmat}
\def\matAABBB[#1,#2,#3;#4,#5,#6]{\beginmat
   #1&#2&#3\\#4&#5&#6\endmat}
\def\matAAAB[#1;#2;#3]{\beginmat#1\\#2\\#3\endmat}
\def\matAAABB[#1,#2;#3,#4;#5,#6]{\beginmat
   #1&#2\\#3&#4\\#5&#6\endmat}
\def\matAAABBB[#1,#2,#3;#4,#5,#6;#7,#8,#9]{\beginmat
   #1&#2&#3\\#4&#5&#6\\#7&#8&#9\endmat}
\def\beginalignorgather#1#2\endalignorgather{%
  \ifx#1!\beginaorgnostar#2\endaorgnostar\else\beginaorgstar#1#2\endaorgstar\fi}
\def\beginaorgstar#1#2\endaorgstar{%
  \ifx#1@\begin{align*}#2\end{align*}\else\begin{gather*}#1#2\end{gather*}\fi}
\def\beginaorgnostar#1#2\endaorgnostar{%
  \ifx#1@\begin{align}#2\end{align}\else\begin{gather}#1#2\end{gather}\fi}
\def\[#1\]{\beginalignorgather#1\endalignorgather}

\def\dbltag#1#2{\tag*{\hbox to 0pt{\hbox to \hsize
  {\hfil#2}\hss}#1}}

\newcommand{\lowsim}{\smash{\hbox{\lower2.5pt
  \hbox{\(\scriptstyle\sim\)}}}}
\def\leq{\leqslant}
\def\geq{\geqslant}
\def\ge{\geq}

\def \kang#1 {{\textcolor{red}{#1}} }
\def \bbR{\mathbb{R}}
\newtheorem{theorem}[subsubsection]{Theorem}
\newtheorem{corollary}[subsubsection]{Corollary}
\newtheorem{lemma}[subsubsection]{Lemma}
\newtheorem{proposition}[subsubsection]{Proposition}

\newtheorem{itheorem}[subsection]{Theorem}
\newtheorem{icorollary}[subsection]{Corollary}

\newtheorem{iconjecture}[subsection]{Conjecture}

\newtheorem*{theorem*}{Theorem}
\newtheorem*{proposition*}{Proposition}
\newtheorem*{conjecture*}{Conjecture}
\theoremstyle{definition}

\newtheorem{definition}[subsubsection]{Definition}
\newtheorem{example}[subsubsection]{Example}

\theoremstyle{remark}
\newtheorem*{remark}{Remark}
\usepackage[all]{xy}
\usepackage{amsmath,amssymb,wasysym}
\makeatother
\newskip\dynsz\newskip\dynsza

\def\DynkinMed{\def\WhiteNodes{\Circle}\def\BlackNodes{\CIRCLE}\dynsz=1ex\dynsza=0.8ex\BWNodes{ooooooooo}}
\def\dynxy{\xymatrix@C=4\dynsz @R=4\dynsz @M=0.6\dynsz}
\def\dynxyn{\xymatrix@C=1.3\dynsz @R=4\dynsz @M=0.6\dynsz}
\def\dynxys{\xymatrix@C=4\dynsz @R=1.0\dynsz @M=0.6\dynsz}
\def\dynxyns{\xymatrix@C=1.3\dynsz @R=1.0\dynsz @M=0.6\dynsz}
\def\xxxcbox#1#2#3{\mbox{\(\vcenter to #1{\vss\vskip 1ex\hbox to #2{\hss#3\hss}\vss}\)}}
\def\xxxubox#1#2#3{\vbox to #1{\vss\hbox to #2{\hss#3\hss}}}
\def\xxxlbox#1#2#3{\vbox to #1{\hbox to #2{\hss#3\hss}\vss}}
\def\myoverunder#1#2#3{\fbox{\(\vcenter to 0.8\dynsz{\vss\baselineskip=0pt\hbox to
    0.8\dynsz{\hss\(\mathop{#1}\limits^{\hbox to
    0ex{\hss\(\vphantom{0}#2\vphantom{0}\)\hss}}_{\hbox to 0pt{\hss\(\vphantom{0}#3\vphantom{0}\)\hss}}\)\hss}\vss}\)}}
\def\myoverunder#1#2#3{\xxxcbox{0.6\dynsz}{0.7\dynsz}{\(\displaystyle\mathop{#1}\limits
^{\xxxubox{0pt}{0pt}{\(#2\)}}%
_{\raise\dynsza\xxxlbox{0pt}{0pt}{\(#3\)}}%
\)}}
\def\dynrefl{\raise0.6\dynsza\hbox{\xxxcbox{0pt}{0pt}{\(\scriptstyle\updownarrow\)}}}
\def\preDynkinNode#1#2#3{\myoverunder{#3}{\scriptstyle#1}{\scriptstyle#2}}
\def\DynkinNode#1{\preDynkinNode{\csname UpperLabel#1\endcsname}%
{\csname LowerLabel#1\endcsname}%
{\csname Node#1\endcsname}%
}
\renewcommand{\xyline}[1][r]{\ar@{-}[#1]}
\renewcommand{\xyto}[1][r]{\ar@{->}[#1]}
\renewcommand{\xytoto}[1][r]{\ar@2{->}[#1]}
\newcommand{\xytototo}[1][r]{\ar@3{->}[#1]}
\newcommand{\xytotototor}{\ar@{=>}[r]\ar@<0.65ex>@{-}[r]+<-1.5ex,0ex>\ar@<-0.65ex>@{-}[r]+<-1.5ex,0ex>}
\newcommand{\xytotototol}{\ar@{=>}[l]\ar@<0.65ex>@{-}[l]+<+1.5ex,0ex>\ar@<-0.65ex>@{-}[l]+<+1.5ex,0ex>}
\newcommand{\xyotot}[1][r]{\ar@2{<-}[#1]}
\newcommand{\xyototot}[1][r]{\ar@3{<-}[#1]}
\def\preSetUpperLabels#1#2#3#4#5#6#7#8#9{%
\def\UpperLabelA{#1}\def\UpperLabelB{#2}\def\UpperLabelC{#3}%
\def\UpperLabelD{#4}\def\UpperLabelE{#5}\def\UpperLabelF{#6}%
\def\UpperLabelG{#7}\def\UpperLabelH{#8}\def\UpperLabelI{#9}}
\def\SetUpperLabels#1{\preSetUpperLabels#1{}{}{}{}{}{}{}{}{}}
\def\preSetLowerLabels#1#2#3#4#5#6#7#8#9{%
\def\LowerLabelA{#1}\def\LowerLabelB{#2}\def\LowerLabelC{#3}%
\def\LowerLabelD{#4}\def\LowerLabelE{#5}\def\LowerLabelF{#6}%
\def\LowerLabelG{#7}\def\LowerLabelH{#8}\def\LowerLabelI{#9}}
\def\SetLowerLabels#1{\preSetLowerLabels#1{}{}{}{}{}{}{}{}{}}
\def\preBWNodes#1#2#3#4#5#6#7#8#9{%
\if#1o\def\NodeA{\WhiteNodes}\else\def\NodeA{\BlackNodes}\fi
\if#2o\def\NodeB{\WhiteNodes}\else\def\NodeB{\BlackNodes}\fi
\if#3o\def\NodeC{\WhiteNodes}\else\def\NodeC{\BlackNodes}\fi
\if#4o\def\NodeD{\WhiteNodes}\else\def\NodeD{\BlackNodes}\fi
\if#5o\def\NodeE{\WhiteNodes}\else\def\NodeE{\BlackNodes}\fi
\if#6o\def\NodeF{\WhiteNodes}\else\def\NodeF{\BlackNodes}\fi
\if#7o\def\NodeG{\WhiteNodes}\else\def\NodeG{\BlackNodes}\fi
\if#8o\def\NodeH{\WhiteNodes}\else\def\NodeH{\BlackNodes}\fi
\if#9o\def\NodeI{\WhiteNodes}\else\def\NodeI{\BlackNodes}\fi}
\def\BWNodes#1{\preBWNodes#1{}{}{}{}{}{}{}{}{}}
\DynkinMed

\newcommand{\DynkinARoAii}[2][]{\SetUpperLabels{#1}\SetLowerLabels{#2}
\mbox{\(\vcenter{\dynxys{
&&\DynkinNode{B}\xyline
&\cdots\xyline
&\DynkinNode{A}\\
\ar@{.}[rrrrr]&\DynkinNode{C}\xyline[ru]\xyline[rd]&&\dynrefl&&\\
&&\DynkinNode{D}\xyline
&\cdots\xyline
&\DynkinNode{E}
}}\)}}
\newcommand{\DynkinAReAii}[2][]{\SetUpperLabels{#1}\SetLowerLabels{#2}
\mbox{\(\vcenter{\dynxys{
&\DynkinNode{C}\xyline\xyline[dd]
&\cdots\xyline
&\DynkinNode{B}\xyline
&\DynkinNode{A}\\
\ar@{.}[rrrrr]&&\dynrefl&&&\\
&\DynkinNode{D}\xyline
&\cdots\xyline
&\DynkinNode{E}\xyline
&\DynkinNode{F}
}}\)}}
\newcommand{\DynkinARiiAii}[2][]{\SetUpperLabels{#1}\SetLowerLabels{#2}
\mbox{\(\vcenter{\dynxy{
\DynkinNode{A}\xyline&\DynkinNode{B}\ar@{<.>}@(u,u)@/_3ex/[l]
}}\)}}
\newcommand{\DynkinARiiiAii}[2][]{\SetUpperLabels{#1}\SetLowerLabels{#2}
\mbox{\(\vcenter{\dynxys@C=2\dynsz{
&&&\DynkinNode{A}\\
\ar@{.}[rrrr]&\DynkinNode{B}\xyline[rru]\xyline[rrd]&&\llap{\dynrefl\;}&&\\
&&&\DynkinNode{C}
}}\)}}
\newcommand{\DynkinARivAii}[2][]{\SetUpperLabels{#1}\SetLowerLabels{#2}
\mbox{\(\vcenter{\dynxyns{
&\DynkinNode{B}\xyline[rr]\xyline[dd]
&&\DynkinNode{A}\\
\ar@{.}[rrrr]&&\dynrefl&&\\
&\DynkinNode{C}\xyline[rr]
&&\DynkinNode{D}
}}\)}}

\newcommand{\DynkinDRnAii}[2][]{\SetUpperLabels{#1}\SetLowerLabels{#2}
\mbox{\(\vcenter{\dynxys{
&&&&&&&\DynkinNode{F}\xyline[ld]\\
\ar@{.}[rrrrrrrr]&\DynkinNode{A}\xyline
&\DynkinNode{B}\xyline
&\DynkinNode{C}\xyline
&\cdots\xyline
&\DynkinNode{D}\xyline
&\DynkinNode{E}&\dynrefl&\\
&&&&&&&\DynkinNode{G}\xyline[lu]
}}\)}}
\newcommand{\DynkinDRivAii}[2][]{\SetUpperLabels{#1}\SetLowerLabels{#2}
\mbox{\(\vcenter{\dynxys{
&&\DynkinNode{C}\xyline[ld]\\
\DynkinNode{A}\xyline&\DynkinNode{B}\\
&&\DynkinNode{D}\xyline[lu]\ar@{<.>}@/_2\dynsz/[uu]
}}\)}}
\newcommand{\DynkinDRivAiii}[2][]{\SetUpperLabels{#1}\SetLowerLabels{#2}
\mbox{\(\vcenter{\dynxys{
&&\DynkinNode{C}\xyline[ld]\ar@{.>}@/_3\dynsz/[lld]\\
\DynkinNode{A}\xyline\ar@{.>}@/_3\dynsz/[rrd]&\DynkinNode{B}\\
&&\DynkinNode{D}\xyline[lu]\ar@{.>}@/_2\dynsz/[uu]
}}\)}}

\newcommand{\DynkinERviAii}[2][]{\SetUpperLabels{#1}\SetLowerLabels{#2}
\mbox{\(\vcenter{\dynxys{
&&&\DynkinNode{C}\xyline&\DynkinNode{A}\\
\ar@{.}[rrrrr]&\DynkinNode{B}\xyline&\DynkinNode{D}\xyline[ru]\xyline[rd]&
\rlap{\hskip 3\dynsz\(\dynrefl\)}&&\\
&&&\DynkinNode{E}\xyline&\DynkinNode{F}
}}\)}}

\newcommand{\DynkinAAffRnAr}[2][]{\SetUpperLabels{#1}\SetLowerLabels{#2}
\mbox{\(\vcenter{\dynxyns@C=0.4\dynsz{
\ar@{.}@/_5.5\dynsz/[dddddddd]&&&&&&&&&&&\DynkinNode{C}\rlap{\;\(\scriptstyle 2\)}\xyline[rrdd]\ar@{.}@/_3\dynsz/[lllllllllll]\\
\\
&&&&&&&&&&&&&\DynkinNode{B}\rlap{\;\(\scriptstyle 1\)}\xyline[rdd]\\
\\
&&&&&&\myoverunder{\circlearrowleft}{}{\genfrac{}{}{0pt}{2}{i \mapsto i+a}{i \in \mathbf{Z}/n\mathbf{Z}}}&&&&&&&&\DynkinNode{A}\rlap{\;\(\scriptstyle 0\)}\\
\\
&&&&&&&&&&&&&\DynkinNode{E}\rlap{\;\(\scriptstyle n-1\)}\xyline[ruu]\\
\\
&&&&&&&&&&&\DynkinNode{D}\rlap{\;\(\scriptstyle n-2\)}\xyline[rruu]\ar@{.}@/^3\dynsz/[lllllllllll]
}}\)}}
\newcommand{\DynkinAAffReAii}[2][]{\SetUpperLabels{#1}\SetLowerLabels{#2}
\mbox{\(\vcenter{\dynxys{
&&\DynkinNode{B}\xyline
&\cdots\xyline
&\DynkinNode{C}\xyline[dd]\\
\ar@{.}[rrrrr]&\DynkinNode{A}\xyline[ru]\xyline[rd]&&\dynrefl&&\\
&&\DynkinNode{E}\xyline
&\cdots\xyline
&\DynkinNode{D}
}}\)}}
\newcommand{\DynkinAAffRoAii}[2][]{\SetUpperLabels{#1}\SetLowerLabels{#2}
\mbox{\(\vcenter{\dynxys{
&&\DynkinNode{B}\xyline
&\cdots\xyline
&\DynkinNode{C}\xyline[rd]\\
\ar@{.}[rrrrrr]&\DynkinNode{A}\xyline[ru]\xyline[rd]
&&\dynrefl&&\DynkinNode{D}\xyline[ld]&\\
&&\DynkinNode{F}\xyline
&\cdots\xyline
&\DynkinNode{E}
}}\)}}
\newcommand{\DynkinAAffRoAiia}[2][]{\SetUpperLabels{#1}\SetLowerLabels{#2}
\mbox{\(\vcenter{\dynxys{
&\DynkinNode{A}\xyline\xyline[dd]
&\DynkinNode{B}\xyline
&\cdots\xyline
&\DynkinNode{C}\xyline
&\DynkinNode{D}\xyline[dd]\\
\ar@{.}[rrrrrr]&&&\dynrefl&&&\\
&\DynkinNode{H}\xyline
&\DynkinNode{G}\xyline
&\cdots\xyline
&\DynkinNode{F}\xyline
&\DynkinNode{E}
}}\)}}
\newcommand{\DynkinAAffRiAii}[2][]{\SetUpperLabels{#1}\SetLowerLabels{#2}
\mbox{\(\vcenter{\dynxy{
\DynkinNode{A}\xyline_\infty&\DynkinNode{B}\ar@{<.>}@(u,u)@/_3ex/[l]
}}\)}}
\newcommand{\DynkinAAffRiiAii}[2][]{\SetUpperLabels{#1}\SetLowerLabels{#2}
\mbox{\(\vcenter{\dynxys@C=2\dynsz{
&&&\DynkinNode{B}\xyline[dd]\\
\ar@{.}[rrrr]&\DynkinNode{A}\xyline[rru]\xyline[rrd]&&\llap{\dynrefl\;}&&\\
&&&\DynkinNode{C}
}}\)}}
\newcommand{\DynkinAAffRiiiAii}[2][]{\SetUpperLabels{#1}\SetLowerLabels{#2}
\mbox{\(\vcenter{\dynxyns{
&&&\DynkinNode{B}\\
\ar@{.}[rrrrrr]&\DynkinNode{A}\xyline[rru]\xyline[rrd]
&&\dynrefl&&\DynkinNode{C}\xyline[lld]\xyline[llu]&\\
&&&\DynkinNode{D}
}}\)}}
\newcommand{\DynkinAAffRiiiAiia}[2][]{\SetUpperLabels{#1}\SetLowerLabels{#2}
\mbox{\(\vcenter{\dynxyns{
&\DynkinNode{A}\xyline[rr]\xyline[dd]
&&\DynkinNode{B}\xyline[dd]\\
\ar@{.}[rrrr]&&\dynrefl&&\\
&\DynkinNode{D}\xyline[rr]
&&\DynkinNode{C}
}}\)}}
\newcommand{\DynkinAAffRivAii}[2][]{\SetUpperLabels{#1}\SetLowerLabels{#2}
\mbox{\(\vcenter{\dynxyns{
&&&\DynkinNode{B}\xyline[rr]
&&\DynkinNode{C}\xyline[dd]\\
\ar@{.}[rrrrrr]&\DynkinNode{A}\xyline[rru]\xyline[rrd]&&&\dynrefl&&\\
&&&\DynkinNode{E}\xyline[rr]
&&\DynkinNode{D}
}}\)}}
\newcommand{\DynkinBAffRnAii}[2][]{\SetUpperLabels{#1}\SetLowerLabels{#2}
\mbox{\(\vcenter{\dynxys{
&\DynkinNode{B}\xyline[rd]\\
\ar@{.}[rrrrrrr]&\dynrefl&\DynkinNode{C}\xyline
&\DynkinNode{D}\xyline
&\cdots\xyline
&\DynkinNode{E}\xytoto
&\DynkinNode{F}&\\
&\DynkinNode{A}\xyline[ur]
}}\)}}

\newcommand{\DynkinBAffRiiiAii}[2][]{\SetUpperLabels{#1}\SetLowerLabels{#2}
\mbox{\(\vcenter{\dynxys{
&\DynkinNode{B}\xyline[rd]\\
\ar@{.}[rrrr]&\dynrefl&\DynkinNode{C}\xytoto
&\DynkinNode{D}&\\
&\DynkinNode{A}\xyline[ur]
}}\)}}
\newcommand{\DynkinBAffRivAii}[2][]{\SetUpperLabels{#1}\SetLowerLabels{#2}
\mbox{\(\vcenter{\dynxys{
&\DynkinNode{B}\xyline[rd]\\
\ar@{.}[rrrrr]&\dynrefl&\DynkinNode{C}\xyline
&\DynkinNode{D}\xytoto
&\DynkinNode{E}&\\
&\DynkinNode{A}\xyline[ur]
}}\)}}
\newcommand{\DynkinBAffDualRnAii}[2][]{
\SetUpperLabels{#1}\SetLowerLabels{#2}
\mbox{\(\vcenter{\dynxys{
&\DynkinNode{B}\xyline[rd]\\
\ar@{.}[rrrrrrr]&\dynrefl&\DynkinNode{C}\xyline
&\DynkinNode{D}\xyline
&\cdots\xyline
&\DynkinNode{E}
&\DynkinNode{F}\xytoto[l]&\\
&\DynkinNode{A}\xyline[ur]
}}\)}}

\newcommand{\DynkinBAffDualRiiiAii}[2][]{\SetUpperLabels{#1}\SetLowerLabels{#2}
\mbox{\(\vcenter{\dynxys{
&\DynkinNode{B}\xyline[rd]\\
\ar@{.}[rrrr]&\dynrefl&\DynkinNode{C}\xyotot
&\DynkinNode{D}&\\
&\DynkinNode{A}\xyline[ur]
}}\)}}
\newcommand{\DynkinBAffDualRivAii}[2][]{\SetUpperLabels{#1}\SetLowerLabels{#2}
\mbox{\(\vcenter{\dynxys{
&\DynkinNode{B}\xyline[rd]\\
\ar@{.}[rrrrr]&\dynrefl&\DynkinNode{C}\xyline
&\DynkinNode{D}\xyotot
&\DynkinNode{E}&\\
&\DynkinNode{A}\xyline[ur]
}}\)}}
\newcommand{\DynkinCAffReAii}[2][]{\SetUpperLabels{#1}\SetLowerLabels{#2}
\mbox{\(\vcenter{\dynxys{
&&\DynkinNode{C}\xyline&\cdots\xyline&\DynkinNode{B}&\DynkinNode{A}\xytoto[l]\\
\ar@{.}[rrrrrr]&\DynkinNode{D}\xyline[ur]\xyline[dr]&&\dynrefl&&&\\
&&\DynkinNode{E}\xyline&\cdots\xyline&\DynkinNode{F}&\DynkinNode{G}\xytoto[l]
}}\)}}
\newcommand{\DynkinCAffRoAii}[2][]{\SetUpperLabels{#1}\SetLowerLabels{#2}
\mbox{\(\vcenter{\dynxys{
&\DynkinNode{C}\xyline\xyline[dd]&\cdots\xyline&\DynkinNode{B}&\DynkinNode{A}\xytoto[l]\\
\ar@{.}[rrrrr]&&\dynrefl&&&\\
&\DynkinNode{D}\xyline&\cdots\xyline&\DynkinNode{E}&\DynkinNode{F}\xytoto[l]
}}\)}}
\newcommand{\DynkinCAffRiiAii}[2][]{\SetUpperLabels{#1}\SetLowerLabels{#2}
\mbox{\(\vcenter{\dynxys@C=2\dynsz{
&&&\DynkinNode{A}\xytoto[lld]\\
\ar@{.}[rrrr]&\DynkinNode{B}&&\dynrefl&\\
&&&\DynkinNode{C}\xytoto[llu]
}}\)}}
\newcommand{\DynkinCAffRiiiAii}[2][]{\SetUpperLabels{#1}\SetLowerLabels{#2}
\mbox{\(\vcenter{\dynxyns{
&\DynkinNode{B}\xyline[dd]&&\DynkinNode{A}\xytoto[ll]\\
\ar@{.}[rrrr]&&\dynrefl&&\\
&\DynkinNode{C}&&\DynkinNode{D}\xytoto[ll]
}}\)}}
\newcommand{\DynkinCAffRivAii}[2][]{\SetUpperLabels{#1}\SetLowerLabels{#2}
\mbox{\(\vcenter{\dynxyns{
&&&\DynkinNode{B}&&\DynkinNode{A}\xytoto[ll]\\
\ar@{.}[rrrrrr]&\DynkinNode{C}\xyline[urr]\xyline[drr]&&&\dynrefl&&\\
&&&\DynkinNode{D}&&\DynkinNode{E}\xytoto[ll]
}}\)}}
\newcommand{\DynkinCAffDualReAii}[2][]{\SetUpperLabels{#1}\SetLowerLabels{#2}
\mbox{\(\vcenter{\dynxys{
&&\DynkinNode{C}\xyline&\cdots\xyline&\DynkinNode{B}\xytoto&\DynkinNode{A}\\
\ar@{.}[rrrrrr]&\DynkinNode{D}\xyline[ur]\xyline[dr]&&\dynrefl&&&\\
&&\DynkinNode{E}\xyline&\cdots\xyline&\DynkinNode{F}\xytoto&\DynkinNode{G}
}}\)}}
\newcommand{\DynkinCAffDualRoAii}[2][]{\SetUpperLabels{#1}\SetLowerLabels{#2}
\mbox{\(\vcenter{\dynxys{
&\DynkinNode{C}\xyline\xyline[dd]&\cdots\xyline&\DynkinNode{B}\xytoto&\DynkinNode{A}\\
\ar@{.}[rrrrr]&&\dynrefl&&&\\
&\DynkinNode{D}\xyline&\cdots\xyline&\DynkinNode{E}\xytoto&\DynkinNode{F}
}}\)}}
\newcommand{\DynkinCAffDualRiiAii}[2][]{\SetUpperLabels{#1}\SetLowerLabels{#2}
\mbox{\(\vcenter{\dynxys@C=2\dynsz{
&&&\DynkinNode{A}\xyotot[lld]\\
\ar@{.}[rrrr]&\DynkinNode{B}&&\dynrefl&\\
&&&\DynkinNode{C}\xyotot[llu]
}}\)}}
\newcommand{\DynkinCAffDualRiiiAii}[2][]{\SetUpperLabels{#1}\SetLowerLabels{#2}
\mbox{\(\vcenter{\dynxyns{
&\DynkinNode{B}\xyline[dd]&&\DynkinNode{A}\xyotot[ll]\\
\ar@{.}[rrrr]&&\dynrefl&&\\
&\DynkinNode{C}&&\DynkinNode{D}\xyotot[ll]
}}\)}}
\newcommand{\DynkinCAffDualRivAii}[2][]{\SetUpperLabels{#1}\SetLowerLabels{#2}
\mbox{\(\vcenter{\dynxyns{
&&&\DynkinNode{B}&&\DynkinNode{A}\xyotot[ll]\\
\ar@{.}[rrrrrr]&\DynkinNode{C}\xyline[urr]\xyline[drr]&&&\dynrefl&&\\
&&&\DynkinNode{D}&&\DynkinNode{E}\xyotot[ll]
}}\)}}
\newcommand{\DynkinDAffRivAii}[2][]{\SetUpperLabels{#1}\SetLowerLabels{#2}
\mbox{\(\vcenter{\dynxys{
&\DynkinNode{B}\xyline[rd]&&\DynkinNode{D}\xyline[ld]\\
\ar@{.}[rrrr]&\dynrefl&\DynkinNode{C}&\dynrefl&\\
&\DynkinNode{A}\xyline[ur]&&\DynkinNode{E}\xyline[lu]
}}\)}}
\newcommand{\DynkinDAffRivAiia}[2][]{\SetUpperLabels{#1}\SetLowerLabels{#2}
\mbox{\(\vcenter{\dynxys@R=0.5\dynsz{
&&\DynkinNode{D}\xyline[ldd]\\
\DynkinNode{B}\xyline[rd]\ar@(ul,dl)@{.>}[]\\
&\DynkinNode{C}\ar@{.}[rr]&\dynrefl&\\
\DynkinNode{A}\xyline[ur]\ar@(ul,dl)@{.>}[]\\
&&\DynkinNode{E}\xyline[luu]
}}\)}}
\newcommand{\DynkinDAffRivAiii}[2][]{\SetUpperLabels{#1}\SetLowerLabels{#2}
\mbox{\(\vcenter{\dynxy{
\DynkinNode{B}\ar@{.>}@/^0pt/[d]\\
\DynkinNode{D}\ar@{.>}@/^0pt/[d]\xyline&\DynkinNode{C}\xyline\xyline[lu]\xyline[ld]&\DynkinNode{A}\\
\DynkinNode{E}\ar@{.>}@/^4\dynsz/[uu]
}}\)}}
\newcommand{\DynkinDAffRivAiv}[2][]{\SetUpperLabels{#1}\SetLowerLabels{#2}
\mbox{\(\vcenter{\dynxy{
\DynkinNode{B}\ar@{.>}@/^1.4\dynsz/[rr]\xyline[rd]&&\DynkinNode{D}\ar@{.>}@/^1.4\dynsz/[dd]\xyline[ld]\\
&\DynkinNode{C}\\
\DynkinNode{A}\ar@{.>}@/^1.4\dynsz/[uu]\xyline[ur]&&\DynkinNode{E}\ar@{.>}@/^1.4\dynsz/[ll]\xyline[lu]
}}\)}}
\newcommand{\DynkinDAffRnAii}[2][]{\SetUpperLabels{#1}\SetLowerLabels{#2}
\mbox{\(\vcenter{\dynxy{
&\DynkinNode{B}\xyline[rd]
&&&&&&\DynkinNode{G}\xyline[ld]\\
\ar@{.}[rr]&\dynrefl&\DynkinNode{C}\xyline
&\DynkinNode{D}\xyline
&\cdots\xyline
&\DynkinNode{E}\xyline
&\DynkinNode{F}&\dynrefl&\ar@{.}[ll]\\
&\DynkinNode{A}\xyline[ur]
&&&&&&\DynkinNode{H}\xyline[lu]
}}\)}}
\newcommand{\DynkinDAffRnAiia}[2][]{\SetUpperLabels{#1}\SetLowerLabels{#2}
\mbox{\(\vcenter{\dynxy@R=0.5\dynsz{
&&&&&&&\DynkinNode{G}\xyline[ldd]\\
&\DynkinNode{B}\xyline[rd]\ar@(ul,dl)@{.>}[]\\
&&\DynkinNode{C}\xyline
&\DynkinNode{D}\xyline
&\cdots\xyline
&\DynkinNode{E}\xyline
&\DynkinNode{F}&\dynrefl&\ar@{.}[ll]\\
&\DynkinNode{A}\xyline[ur]\ar@(ul,dl)@{.>}[]\\
&&&&&&&\DynkinNode{H}\xyline[luu]
}}\)}}
\newcommand{\DynkinDAffReAii}[2][]{\SetUpperLabels{#1}\SetLowerLabels{#2}
\mbox{\(\vcenter{\dynxys{
&&&&&\DynkinNode{A}\ar@{<.>}@/^4\dynsz/[dddddd]\\
&&\DynkinNode{D}\xyline&\cdots\xyline&\DynkinNode{C}\xyline[ru]\xyline[rd]\\
&&&&&\DynkinNode{B}\ar@{<.>}@/^1.5\dynsz/[dd]\\
\ar@{.}[rrrrrr]&\DynkinNode{E}\xyline[ruu]\xyline[rdd]&&\dynrefl&&&\\
&&&&&\DynkinNode{H}\\
&&\DynkinNode{F}\xyline&\cdots\xyline&\DynkinNode{G}\xyline[ru]\xyline[rd]\\
&&&&&\DynkinNode{I}
}}\)}}
\newcommand{\DynkinDAffRoAii}[2][]{\SetUpperLabels{#1}\SetLowerLabels{#2}
\mbox{\(\vcenter{\dynxys{
&&&&\DynkinNode{A}\ar@{<.>}@/^4\dynsz/[dddddd]\\
&\DynkinNode{D}\xyline\xyline[dddd]&\cdots\xyline&\DynkinNode{C}\xyline[ru]\xyline[rd]\\
&&&&\DynkinNode{B}\ar@{<.>}@/^1.5\dynsz/[dd]\\
\ar@{.}[rrrrr]&&\dynrefl&&&\\
&&&&\DynkinNode{G}\\
&\DynkinNode{E}\xyline&\cdots\xyline&\DynkinNode{F}\xyline[ru]\xyline[rd]\\
&&&&\DynkinNode{H}
}}\)}}
\newcommand{\DynkinDAffReAiv}[2][]{\SetUpperLabels{#1}\SetLowerLabels{#2}
\mbox{\(\vcenter{\dynxys{
&&&&&\DynkinNode{A}\ar@{.>}@/_4\dynsz/[dddddd]\\
&&\DynkinNode{D}\xyline&\cdots\xyline&\DynkinNode{C}\xyline[ru]\xyline[rd]\\
&&&&&\DynkinNode{B}\ar@{.>}@/_1.5\dynsz/[dd]\\
\ar@{.}[rrrrrr]&\DynkinNode{E}\xyline[ruu]\xyline[rdd]&&\dynrefl&&&\\
&&&&&\DynkinNode{H}\ar@{.>}@/_7\dynsz/[uuuu]\\
&&\DynkinNode{F}\xyline&\cdots\xyline&\DynkinNode{G}\xyline[ru]\xyline[rd]\\
&&&&&\DynkinNode{I}\ar@{.>}@/_7\dynsz/[uuuu]
}}\)}}
\newcommand{\DynkinDAffRoAiv}[2][]{\SetUpperLabels{#1}\SetLowerLabels{#2}
\mbox{\(\vcenter{\dynxys{
&&&&\DynkinNode{A}\ar@{.>}@/_4\dynsz/[dddddd]\\
&\DynkinNode{D}\xyline\xyline[dddd]&\cdots\xyline&\DynkinNode{C}\xyline[ru]\xyline[rd]\\
&&&&\DynkinNode{B}\ar@{.>}@/_1.5\dynsz/[dd]\\
\ar@{.}[rrrrr]&&\dynrefl&&&\\
&&&&\DynkinNode{G}\ar@{.>}@/_7\dynsz/[uuuu]\\
&\DynkinNode{E}\xyline&\cdots\xyline&\DynkinNode{F}\xyline[ru]\xyline[rd]\\
&&&&\DynkinNode{H}\ar@{.>}@/_7\dynsz/[uuuu]
}}\)}}
\newcommand{\DynkinEAffRviAii}[2][]{\SetUpperLabels{#1}\SetLowerLabels{#2}
\mbox{\(\vcenter{\dynxys{
&&&&\DynkinNode{D}\xyline&\DynkinNode{B}\\
\ar@{.}[rrrrrr]&\DynkinNode{A}\xyline&\DynkinNode{C}\xyline&\DynkinNode{E}\xyline[ru]\xyline[rd]&
\rlap{\hskip 3\dynsz\(\dynrefl\)}&&\\
&&&&\DynkinNode{F}\xyline&\DynkinNode{G}
}}\)}}
\newcommand{\DynkinEAffRviAiii}[2][]{\SetUpperLabels{#1}\SetLowerLabels{#2}
\mbox{\(\vcenter{\dynxys{\\
\DynkinNode{B}\ar@{.>}@/_6\dynsz/[rrdddddd]&&&&\DynkinNode{G}\ar@{.>}@/_6\dynsz/[llll]\\
&\DynkinNode{D}\xyline[lu]\ar@{.>}@/_3\dynsz/[rddd]&&\DynkinNode{F}\xyline[ru]\ar@{.>}@/_3\dynsz/[ll]\\
&&\DynkinNode{E}\xyline[lu]\xyline[ru]\\ \\
&&\DynkinNode{C}\xyline[uu]\ar@{.>}@/_3\dynsz/[ruuu]\\ \\
&&\DynkinNode{A}\xyline[uu]\ar@{.>}@/_6\dynsz/[rruuuuuu]
}}\)}}
\newcommand{\DynkinEAffRviiAii}[2][]{\SetUpperLabels{#1}\SetLowerLabels{#2}
\mbox{\(\vcenter{\dynxys{
&&\DynkinNode{D}\xyline&\DynkinNode{B}\xyline&\DynkinNode{A}\\
\DynkinNode{C}\xyline&\DynkinNode{E}\xyline[ur]\xyline[dr]\ar@{.}[rrrr]&&\dynrefl&&\\
&&\DynkinNode{F}\xyline&\DynkinNode{G}\xyline&\DynkinNode{H}
}}\)}}
\def\starwedge#1#2#3#4|{#1\if#2*Aff\fi\if#3^Dual\fi}
\def\rankauto#1[#2]#3|{Rank\if#1NN\else\romannumeral#1\fi Auto\ifnum#2<0 iiA\else\romannumeral#2\fi}
\def\Dynkin#1#2{\edef\aaa{\starwedge#1--|\rankauto#2[1]-|}\aaa\csname\aaa\endcsname}

\def\starwedge#1#2#3#4#5|{%
  \if#1B%
    \if#2C#1#2\if#3*Aff\fi\if#3^Dual\fi\if#4^Dual\fi
    \else#1\if#2*Aff\fi\if#2^Dual\fi\if#3^Dual\fi
    \fi
  \else
    #1\if#2*Aff\fi\if#2^Dual\fi\if#3^Dual\fi
  \fi}
\def\rankauto#1[#2]#3|{R\if#1nn\else\if#1oo\else\if#1ee\else\romannumeral#1\fi\fi\fi 
A\if#2rr\else\ifnum#2<0 iia\else\romannumeral#2\fi\fi}
\def\Dynkin#1#2{\edef\dynkincsname{Dynkin\starwedge#1----|\rankauto#2[1]-|}%
\csname\dynkincsname\endcsname}
\def\removearrows{\renewcommand{\xyto}[1][r]{\ar@{-}[##1]}%
\renewcommand{\xytoto}[1][r]{\ar@2{-}[##1]}%
\renewcommand{\xytototo}[1][r]{\ar@3{-}[##1]}
\renewcommand{\xytotototor}{\ar@{=}[r]\ar@<0.65ex>@{-}[r]\ar@<-0.65ex>@{-}[r]}%
\renewcommand{\xytotototol}{\ar@{=}[l]\ar@<0.65ex>@{-}[l]\ar@<-0.65ex>@{-}[l]}%
\renewcommand{\xyotot}[1][r]{\ar@2{-}[##1]}%
\renewcommand{\xyototot}[1][r]{\ar@3{-}[##1]}}
\def\restorearrows{\renewcommand{\xyto}[1][r]{\ar@{->}[##1]}%
\renewcommand{\xytoto}[1][r]{\ar@2{->}[##1]}%
\renewcommand{\xytototo}[1][r]{\ar@3{->}[##1]}%
\renewcommand{\xytotototor}{\ar@{=>}[r]\ar@<0.65ex>@{-}[r]+<-1.5ex,0ex>\ar@<-0.65ex>@{-}[r]+<-1.5ex,0ex>}%
\renewcommand{\xytotototol}{\ar@{=>}[l]\ar@<0.65ex>@{-}[l]+<+1.5ex,0ex>\ar@<-0.65ex>@{-}[l]+<+1.5ex,0ex>}%
\renewcommand{\xyotot}[1][r]{\ar@2{<-}[##1]}%
\renewcommand{\xyototot}[1][r]{\ar@3{<-}[##1]}}

\makeop{BRD}

\def\aff{{\rm aff}}

\def\Waff{W}

\normalsize
\def \Stab {\mathrm{Stab}}
\def \ga {\gamma}

\def \scrC {\mathcal{C}}
\def \scrS {\mathcal{S}}
\def \scrF {\mathcal{F}}
\def \scrP {\mathcal{P}}
\def \scrT {\mathcal{T}}
\def \scrB {\mathcal{B}}
\def \scrA {\mathcal{A}}
\def \scrL {\mathcal{L}}
\def \scrD {\mathcal{D}}

\def \scrR {\mathcal{R}}
\def \scrW {\mathcal{W}}

\title[Hecke Algebra-valued Series and Affine Weyl Groups]{Hecke Algebra-valued Poincar\'e Series and Geometric Factorization of Affine Weyl Groups}
\author{Ming-Hsuan Kang}
\thanks{Email: mhkang@math.nctu.edu.tw}

\author{Jiu Kang Yu}
\thanks{Email: jkyu@ims.cuhk.edu.hk}
\date{}
\begin{document}
\maketitle 

\begin{abstract}
This paper studies affine Weyl groups through the geometry of geodesic
tubes and their hyperbolic stabilizers. We introduce a geometric
framework relating affine Weyl group factorizations to
Hecke-algebra-valued Poincar\'e series and zeta functions on quotient
complexes. Using this framework, we prove a conjecture describing the
connection between geodesic tubes and factorizations of Poincar\'e
series for affine Weyl groups. We also establish partial results toward
a Hecke-algebra-valued zeta identity for simply connected groups,
including the cases of types $\widetilde{A}_{n-1}$ and
$\widetilde{C}_n$. These results highlight the interaction between
affine Weyl group geometry, Hecke algebras, and zeta functions arising
from Bruhat--Tits buildings.
\end{abstract}

\section{Introduction}
\def \h {\text{ht}}

Poincar\'e series of Coxeter groups, and their Hecke-algebra-valued
analogues, provide a natural meeting point for algebraic, geometric,
and combinatorial structures.  In the affine case, these series also
interact with the geometry of affine Coxeter complexes and
Bruhat-Tits buildings.  The purpose of this paper is to make this
interaction explicit through the geometry of geodesic tubes and their
hyperbolic stabilizers.

Let \((W,S)\) be a Coxeter group, and let \(H\) be its Hecke algebra
over \(\bbC\) with parameter \(q \in \bbC^\times\).  We write
\(\{e_w\mid w\in W\}\) for the standard basis of \(H\).  For a subset
\(Y\) of \(W\), there are two Poincar\'e series naturally associated
with \(Y\):
\[
    p_Y := \sum_{x \in Y} u^{\ell(x)} \in \mathbb{Z}[[u]]  \quad \text{and} \quad P_Y := \sum_{w \in Y} e_w u^{\ell(w)} \in H[[u]]
\]
The scalar series \(p_Y\) is classical and has many applications, for
example in invariant theory.  The Hecke-algebra-valued series \(P_Y\)
retains the Coxeter length information together with the Hecke
operators \(e_w\), and \(p_Y\) is obtained from it by specialization.
The central point of the paper is that, for affine Weyl groups,
identities among such series reflect geometric factorizations coming
from hyperbolic elements stabilizing geodesic tubes.

We begin with the familiar \(1\)-dimensional model.  Let \(X\) be a
finite graph obtained as a quotient of the Bruhat-Tits tree of
\(\mathrm{SL}_2(F)\), where \(F\) is a non-archimedean local field.
The \textit{Ihara zeta function} of \(X\) is
\[
Z_X(u) := \exp\left(\sum_{m \geq 1} \frac{N_m(X)}{m} u^m\right) \in \mathbb{Q}[[u]],
\]
where \(N_m(X)\) is the number of geodesic cycles of length \(m\) in
\(X\).  The affine Weyl group of \(\mathrm{SL}_2(F)\) is of type
\(\tilde{A}_1\), with \(S=\{s_1,s_2\}\).  The space
\(\mathbb{C}(\mathrm{Ch}(X))\) of functions on the edges of \(X\)
is naturally a module for the corresponding Hecke algebra \(H\), for a
suitable \(q\).  Hashimoto \cite{Ha1} showed that
\[
Z_X(u) = \det\left(I- e_{s_1s_2} u^2 \mid \mathbb{C}(\mathrm{Ch}(X)) \right)^{-1}.
\]
Hoffman \cite{Hof} observed that this determinant is also expressible
as an alternating product of Hecke-algebra-valued Poincar\'e series over parabolic
subgroups; all operators below act on \(\bbC({\rm Ch}(X))\):
\[
\det(I- e_{s_1s_2} u^2)^{-1} = \det (P_{W} ) \det(P_{W_{\{s_1\}}}
)^{-1} \det(P_{W_{\{s_2\}}} )^{-1} .
\]
The group-theoretic source of this identity is the length-preserving
factorization
\[
W = W_{\{s_2\}}\cdot \{1, s_1s_2, (s_1s_2)^2, \dots\}\cdot W_{\{s_1\}}.
\]
Hoffman asked whether the resulting relation
\[
Z_X(u) = \prod_{I \subset S} \det\big(P_{W_I} \mid  \mathbb{C}(\mathrm{Ch}(X)) \big)^{(-1)^{|S\setminus I|}},
\]
might extend to simply connected \(p\)-adic groups of higher rank.
The main difficulty is that, in higher-dimensional Bruhat-Tits
buildings, one must identify the appropriate higher-dimensional
replacement for geodesic cycles and relate it to affine Weyl group
geometry.

The rank two case was studied in \cite{KM}.  Let \(G\) be a split,
simple, simply connected algebraic group over a non-archimedean local
field \(F\) of rank \(2\).
Let $\omega_1^\vee$ and $\omega_2^\vee$ be the two fundamental
coweights.  In \cite{KM}, geodesic strip zeta functions of direction
\(\omega_i^\vee\), denoted by \(Z^{(\omega_i^\vee)}_X(u)\), are
defined for a finite quotient complex \(X\) of the building of \(G\).
These are the rank-two instances of the geodesic tube zeta functions
considered below.
These zeta functions satisfy an analogue of Hashimoto's formula:
\[
Z^{(\omega_i^\vee)}_X(u)=\det\big(I- e_{w_i}u^{\ell(w_i)}\mid \mathbb{C}(\mathrm{Ch}(X)) \big)^{-1}.
\]
Here \(w_i\) is the straight generator of the corresponding geodesic
strip of type \(\omega_i^\vee\).  The product of the two zeta functions is then
related to the alternating product of Poincar\'e series by
\[
\prod_{i=1}^2 Z^{(\omega_i^\vee)}_X(u) = \prod_{I \subset S} \det (P_{W_I}\mid  \mathbb{C}(\mathrm{Ch}(X)))^{(-1)^{|S\setminus I|}},
\]
through a length-preserving factorization of the affine Weyl group.

This suggests a higher-rank picture in which each fundamental
coweight direction should give rise to a suitable geodesic tube, and
the corresponding hyperbolic stabilizer should supply one factor in
the Poincar\'e-series identity.  The numerical evidence comes from the
scalar Poincar\'e series.  The authors of \cite{KM} showed that there
exist positive integers \(d_1,\ldots,d_n\) such that
\[
\prod_{I \subset S} p_{W_I}\!^{(-1)^{|S\setminus I|}} = \frac{1}{(1-u^{d_1})\cdots (1-u^{d_{n}})}.
\]
These integers, listed in Table 1, predict the lengths of the
straight generators that should occur in the geometric
theory.
\begin{center}
 \begin{tabular}{|c|c|l|}
 \hline
 Type of $W^f$ & Coxeter number $h$ & $d_1, \cdots, d_{n}$ \\ [0.5ex]
 \hline \hline
 ${A_{n}}$ & $n+1$ & $n+1,\cdots,n+1$ \\ \hline
 $B_{n},C_{n}$ & $2n$ & $n+1,n+2,\cdots,2n$ \\ \hline
 ${D_{n}}$ & $2n-2$ & $n+1,n+2, \cdots, 2n-2, 2n-2,2n-2$ \\ \hline
 ${E_{6}}$ & 12 & $7, 9, 9, 11, 12, 12$ \\ \hline
 ${E_{7}}$ & 18 & $8, 10, 11, 13, 14, 17, 18$ \\ \hline
 ${E_{8}}$ & 30 & $9, 11, 13, 14, 17, 19, 23, 29$ \\ \hline
 ${F_{4}}$ & 12 & $5, 7, 8, 11$ \\ \hline
 ${G_{2}}$ & 6 & $3, 5$ \\ \hline
 \end{tabular}\\
 \noindent{\smallbreak}
 Table 1
\end{center}

Motivated by this prediction, \cite{KM} formulated the following
conjectures.  The first is geometric: the integers in Table 1 should
be realized as the lengths of straight generators of geodesic tubes.

\begin{iconjecture} \label{conj-a}
There exist \( w_1,\cdots, w_{n} \in W \) such that \( w_i \) is the straight generator of some geodesic tube and has length \( d_i \). Then
\[
\prod_{I \subset S} p_{W_I}\!^{(-1)^{|S\setminus I|}} = \prod_{i=1}^{n} p_{\bbT_i} \qquad \text{in} \qquad \bbC[[u]]^\times,
\]
where  \( \bbT_i = \{w_i^N: N \geq 0\} \) is the monoid generated by \( w_i \).
\end{iconjecture}

The second conjecture is the Hecke-algebra-valued analogue of the
same factorization.

\begin{iconjecture} \label{conj-b}
For the monoid \( \bbT_i\) in Conjecture \ref{conj-a}, we have 
\[
\prod_{I \subset S} P_{W_I}\!^{(-1)^{|S\setminus I|}} = \prod_{i=1}^{n} P_{\bbT_i} \qquad \text{in} \qquad \left( H[[u]]^\times \right)^{\text{ab}}.
\]
Here $\left(H[[u]]^\times \right)^{\text{ab}}$ is the maximal abelian quotient of $H[[u]]^\times$.
\end{iconjecture}

One point deserves emphasis.  In higher rank, the notion of a
geodesic tube, and especially the structure of its hyperbolic
stabilizer in the affine Weyl group, is part of the structure to be
explained rather than background terminology.  The first contribution of this paper is
to develop this geometry systematically.  We then use it to prove
Conjecture \ref{conj-a}, as stated in Theorem \ref{main1}, and to
prove Conjecture \ref{conj-b} in the cases described in Theorem
\ref{main2}.

\begin{itheorem}\label{main1}
Let \( W \) be an affine Weyl group associated with an irreducible root
system \( \Phi \) of rank \(n\).  Let \(B\) be a system of simple roots of
\(\Phi\).  For any \(\beta\in B\), put
\[
d_\beta = \frac{(\omega_\beta^\vee, \beta^\vee)}{(\omega_\beta^\vee,\omega_\beta^\vee)}\langle 2\rho, \omega_\beta^\vee \rangle,
\]
where \(\beta^\vee\) is the coroot of \(\beta\), \(\omega_\beta^\vee\)
is the fundamental coweight of \(\beta\), and \(\rho\) is the half-sum
of positive roots.  Then \(\{d_\beta\}_{\beta\in B}=\{d_1,\ldots,d_{n}\}\).
Moreover, \(d_\beta\) is the length of the straight
generator $w_\beta$ of a suitable geodesic tube of direction \(
\omega_\beta^\vee \), and we have 
\begin{equation}\label{eqmain1}
\prod_{I \subset S} p_{W_I}\!^{(-1)^{|S\setminus I|}} =
\prod_{\beta\in B}p_{\bbT_\beta}, \quad\mbox{ where }
\bbT_\beta=\{w_\beta^N:N\geq 0\}.
\end{equation}
\end{itheorem}

\begin{itheorem}\label{main2}
Conjecture \ref{conj-b} is valid when \( W \) is either of type \( \tilde{A}_{n-1} \) or \( \tilde{C}_{n} \).
\end{itheorem}

The same geometry also supplies the higher-dimensional cycles needed
for the zeta side.  On a finite quotient complex \(X\), we define
circular geodesic tubes of type \(\beta\), of direction
\(\omega_\beta^\vee\), in analogy with the geodesic strips of
\cite{KM}.  The associated zeta function
\(Z^{(\beta)}_X(u)\) is defined in Section 4.  Its
determinant expression is the following.
\begin{itheorem}\label{main3}
Assuming \( X \) is finite, for each \(\beta \in B \), we have:
\[
Z^{(\beta)}_X(u) = \det(P_{\bbT_\beta} \mid \bbC({\rm Ch}(X))).
\]
\end{itheorem}

Combining Theorem \ref{main2} with Theorem \ref{main3} gives the
following zeta identity for the simply connected case.

\begin{icorollary}[Zeta Identity for Simply Connected Groups]\label{main4}
Assuming \( X \) is finite and $W$ is of type $\tilde{A}_{n-1}$ or $\tilde{C}_{n}$, the following relation holds:
\[
\prod_{\beta\in B} Z^{(\beta)}_X(u) = \prod_{I \subset S} \det (P_{W_I}\mid  \bbC(\mathrm{Ch}(X)))^{(-1)^{|S\setminus I|}}.
\]
\end{icorollary}

In this corollary, \(X\) is a quotient of the building of the simply
connected algebraic group \(G\), and the zeta functions concern
top-dimensional simplices.  This differs from the zeta identities
studied for complexes attached to adjoint type groups.  For such
results, see \cite{KL1} and \cite{KL2} for PGL$_3$, \cite{FLW} for
PGSP$_4$, \cite{DK} for PGL$_n$ over a 1-adic field, and \cite{KLW}
for rank two algebraic groups over a 1-adic field.  In the adjoint
case, the zeta identity expresses an alternating product of zeta
functions in various dimensions as a single unramified Langlands
\(L\)-function.

The paper is organized as follows.  Section 2 develops the
Hecke-algebra-valued Poincar\'e series \(P_Y\), including the
length-preserving decompositions and abelian quotient identities used
later.  Section 3 develops the theory of geodesic tubes in affine
Weyl group geometry, analyzes their hyperbolic stabilizers, and proves
Theorem \ref{main1}.  Section 4 defines the zeta function of circular
geodesic tubes and proves Theorem \ref{main3}.  Sections 5 and 6
prove Theorem \ref{main2} by giving explicit length-preserving
factorizations for affine Weyl groups of types \(\tilde{A}_{n-1}\) and
\(\tilde{C}_{n}\), respectively.

The authors acknowledge the assistance of ChatGPT in improving the
language and presentation of the manuscript.

\section{Poincar\'e Series with Values in the Hecke Algebra}
This section develops the algebraic formalism underlying the geometric factorizations used later.

The Poincar\'e series of (a subset of) a Coxeter system \((W,S)\)
\cite{Mac} is a well-known and useful notion.  It is a power
series with integer coefficients.  We define a Hecke-algebra-valued
variant for \((W,S)\).  The
Poincar\'e series with matrix coefficients introduced by \cite{Gyoja}
and used by \cite{KM} are incarnations of our variant.

\subsection{Definitions}

Let \((W,S)\) be a Coxeter system with \(S\) finite,  and let \(H\) be its
Hecke algebra over a commutative ring \(A\supset \bbQ\) relative to the parameter \(q \in
A^\times\).  By \cite[IV.2, Exercise 23]{Bou}, \(H\) has an
\(A\)-basis \(\{e_w\}_{w\in W}\) and each \(e_w\) is invertible in
\(H\).

For any subset \(X\) of \(W\), the corresponding Poincar\'e series is defined as
\[
  P_X := \sum_{x \in X} e_x u^{\ell(x)} \in H[[u]] \subset H(\!(u)\!) \quad \mbox{and} \quad
  p_X := \sum_{x \in X} u^{\ell(x)} \in \bbZ[[u]] \subset \bbZ(\!(u)\!).
\]
More generally, for an element \(M = \{m_x\}_{x \in W} \in \bbZ[[W]] := \text{Map}(W, \bbZ)\), the set of integer-valued functions on \(W\), we define
\[
  P_M := \sum_{x \in W} m_x e_x u^{\ell(x)} \in H[[u]] \quad \mbox{and} \quad 
  p_M := \sum_{x \in W} m_x u^{\ell(x)} \in \bbZ[[u]].
\]
If \(m_x \in\{0,1\}\) for all \(x \in W\), \(P_M\) reduces to
\(P_X\), where \(X=\{x\in W: m_x=1\}\).  In this case, we often identify
\(X\) and \(M\), and call \(M\) a \textit{set}.  If \(m_x \geq 0\) for all \(x \in W\),
\(M\) is a \textit{multiset} over \(W\).

 Finally, note that \((\rho_1)_*(P_M(u)) = p_M(qu)\), where \(\rho_1: H \to A\) is the ring homomorphism that maps each \(e_s\) to \(q\) for all \(s \in S\), and \((\rho_1)_*\) denotes coefficientwise application of \(\rho_1\).

\subsection{Length-Preserving Decompositions}

We next record the compatibility between multiplication of multisets
and multiplication of their Poincar\'e series.

Let \(M = \{m_x\}_{x \in W}\) and \(M' = \{m'_y\}_{y \in W}\) be two multisets defined over \(W\). Assume that the expression
\[
m''_z := \sum_{\substack{x, y \in W \\ xy = z}} m_x m'_y
\]
has finitely many non-zero terms, so \(m''_z\) is well-defined for each \(z \in W\). The resulting multiset \(M'' = \{m''_z\}_{z \in W}\) is the \textit{product} of \(M\) and \(M'\), denoted \(M'' = M \times M'\). We also call this a \textit{decomposition} of \(M''\).

A decomposition \(M'' = M \times M'\) is said to be \textit{length-preserving} if \(\ell(xy) = \ell(x) + \ell(y)\) holds for all \(x, y \in W\) with \(m_x > 0\) and \(m'_y > 0\).

\begin{example}[\cite{Hum}, \S 5.12]
Let \( I \subseteq J \subseteq S \). There is a unique subset of \( W_J \), denoted by \( W_J/W_I \) (respectively, \( W_I \backslash W_J \)), such that \( (W_J/W_I) \times W_I \) (respectively, \( W_I \times (W_I \backslash W_J) \)) is a length-preserving decomposition of \( W_J \).
\end{example}

\begin{example}\label{exampleAn}
Let \( (W, S) \) be a Coxeter system of type \(A_{n-1}\), and let \( S
= \{ s_1, \ldots, s_{n-1} \} \). Consider \( I = \{ s_2, \ldots,
s_{n-1} \} \), so that \( W \simeq S_{n}, W_I\simeq S_{n-1} \). Then, we have
\[
W/W_I = \{ 1, s_1, s_2s_1, \ldots, s_{n-1}s_{n-2}\cdots s_1 \}\]
and 
\[W_I \backslash W = \{ 1, s_1, s_1s_2, \ldots, s_1s_2\cdots s_{n-1} \}.
\]
\end{example}

\begin{proposition} \label{lengthprevseringdecompsition}
Assume that \(M, M'\) are multisets over \(W\) and \(M'' = M \times M'\) is well-defined. Then the following conditions are equivalent:
\begin{itemize}
\item[\textnormal{(a)}] \(p_{M''} = p_M p_{M'}\) in \(\bbZ(\!(u)\!)\).
\item[\textnormal{(b)}] \(P_{M''} = P_M P_{M'}\) in \(H(\!(u)\!)\).
\item[\textnormal{(c)}] The decomposition \(M'' = M \times M'\) is length-preserving.
\end{itemize}
\end{proposition}

\begin{proof} 
The implication (c) \(\Rightarrow\) (b) follows from the definition, and (b) \(\Rightarrow\) (a) follows from \((\rho_1)_*(P_M(u)) = p_M(qu)\). It remains to prove (a) \(\Rightarrow\) (c).

 Assume (a).  We claim: for any integer \(n\), for any \(x,y \in W\)
 such that \(m_x > 0\), \(m'_y > 0\)
 \[
  \ell(x)+\ell(y)=n \Rightarrow  \ell(xy)=n \qquad\mbox{ and }\qquad \ell(x)+\ell(y)
 > n \Rightarrow \ell(xy) > n.
\]
We prove the claim by induction on \(n\).  The claim is immediate
when \(n<0\).  Suppose  \(n \geq 0\).  Assume
 \(\ell(x)+\ell(y) = n\).  Then \(\ell(xy)=n\) since \(\ell(xy) > n-1\) by the
 induction hypothesis.  Now assume \(m_{x_0} > 0\), \(m'_{y_0}>0\), \(\ell(x_0)+\ell(y_0) > n\) but \(\ell(x_0y_0)\leq
 n\).  Then \(\ell(x_0y_0)=n\) by induction hypothesis.   Now
\[
\sum _{\ell(z)=n} m''_z \geq m_{x_0}m'_{y_0}+\sum_{\ell(x)+\ell(y)=n} m_x m'_y>\sum_{\ell(x)+\ell(y)=n} m_x m'_y.
\]
That is, the coefficient of \(u^n\) in \(p_{M''}\) is greater than that in
\(p_M p_{M'}\), contradicting assumption (a).  This completes the proof
of the claim and the proposition.
\end{proof}

\begin{remark}$ $
\begin{enumerate}
\item[(a)] The multiplication of multisets is associative, i.e., \((M
  \times M') \times M'' = M \times (M' \times M'')\), provided that these multiplications are well-defined.
\item[(b)] \(M \times M'\) is well-defined if \(M\) is finite, meaning
  \(m_x \neq 0\) for only finitely many \(x \in W\).  It is also
  well-defined when \(M'\) is finite.
\end{enumerate}
\end{remark}

\begin{proposition}\label{uniq-dec}
   Let \(L\) be a subset of \(W\) containing the identity element \(e\).
   Assume that \(M,M'\) are finite subsets of \(W\) such that
   \(M\times L=M'\times L\) and the decompositions \(M\times L\),
   \(M'\times L\) are length-preserving.  Then \(M=M'\).
\end{proposition}

\begin{proof}  Replacing \(M\) (resp.~\(M'\)) with \(M\setminus M\cap
  M'\) (resp.~\(M'\setminus M\cap M'\)), we may assume that \(M\) and
  \(M'\) are disjoint.  Assume that they are non-empty and take \(x_0\in
  M\).  Then \(x_0e=x_0'y\) for some \(x_0' \in M'\), \(y
  \in L\).  Then \(l(x_0)=\ell(x_0')+\ell(y)>\ell(x_0')\).  Similarly,
  \(x_0'e=x_1z\) for some \(x_1 \in M\), \(z\in L\), and \(\ell(x_0')>\ell(x_1)\).
  Inductively we have an infinite sequence \(\{x_i\}_{i\geq 0}\) in
  \(M\) with \(\{\ell(x_i)\}_{i\geq0}\) strictly decreasing. This
  contradiction proves the proposition.
\end{proof}

\subsection{Identities in Abelian Quotients of \(H[[u]]^\times\) and \(H(\!(u)\!)^\times\)}

The length-preserving identities obtained above will later be used
after passing to abelian quotients, where products may be reordered.

For a given group \(G\), we denote its maximal abelian quotient by \(G^{\text{ab}}\), which is the quotient of \(G\) by its commutator subgroup.

Consider a multiset \(M = \{m_x\}_{x \in W} \in \bbZ[[W]]\) with a unique element of minimal length. Thus, for some \(x_0\), we have \(m_{x_0} > 0\) and \(\ell(x) > \ell(x_0)\) for all \(x \neq x_0\) with \(m_x > 0\). Then \(P_M\) is invertible in \(H(\!(u)\!)\).

We consider identities of the form
\[
P_{M_1} \cdots P_{M_s} = P_{M_1'} \cdots P_{M_t'} \quad \text{in} \, \left(H(\!(u)\!)^\times\right)^{\text{ab}},
\]
where \(M_1, \ldots, M_s\) and \(M_1', \ldots, M_t'\) each have a unique element of minimal length. Such an identity induces identities in \(A(\!(u)\!)^\times\): let \(V\) be an \(H\)-module whose underlying \(A\)-module is free of finite rank, and let \(\rho:H \to \End_A(V)\) be the associated ring homomorphism. Then
\[
\det(\rho_*P_{M_1}) \cdots \det(\rho_*P_{M_s}) = \det(\rho_*P_{M_1'}) \cdots \det(\rho_*P_{M_t'}) \quad \text{in} \quad A(\!(u)\!)^\times.
\]

If all the multisets \(M_i\) and \(M'_j\) contain the identity element of \(W\), then \(P_{M_i}\) and \(P_{M'_j}\) are invertible in \(H[[u]]\). The identity \(P_{M_1} \cdots P_{M_s} = P_{M_1'} \cdots P_{M_t'}\) may then be considered in \(\left(H[[u]]^\times\right)^{\text{ab}}\), giving identities in \(A[[u]]^\times\).

Retain the notations of
Example~\ref{exampleAn}.  The following relations in
\((H[[u]]^\times)^{\rm ab}\) of the above kind will be used in the
proof of the main result of \ref{conjAn}.

For \(k=0,\ldots,n\), we put \(Y_k=\{\sigma
  _{k,0},\ldots,\sigma_{k,n-1}\}\), where
\[
  \sigma _{k,j}=\begin{cases}
    s_j s_{j-1}\cdots s_1&\mbox{ if }j\leq k,\\
    (s_ks_{k-1}\cdots s_1)(s_{k+1}s_{k+2}\cdots s_j) &\mbox{ if }j\geq
    k.
  \end{cases}
\]
Notice that \(Y_0=Y_1=W_I\backslash W\), \(Y_{n-1}=Y_n=W/W_I\).

\begin{theorem} \label{idenAn} With the above notation, we have
\[
  P_{Y_0}=\cdots=P_{Y_n}\quad\mbox{ in }\quad (H[[u]]^\times)^{\rm
   ab},
\]
and
\[
  P_{W}=P_{W_I}P_{Y_k}\quad\mbox{ in }\quad (H[[u]]^\times)^{\rm
    ab}\qquad \mbox{ for }k=0,\ldots,n.
\]
\end{theorem}

\begin{proof}  By Example~\ref{exampleAn}, we have \(P_W=P_{W_I}P_{Y_k}\) in
  \((H[[u]]^\times)^{\rm ab}\) for \(k=0,1,n-1,n\).
Now assume \(2\leq k\leq n-2\).  Write \(Y_k=Y_k'\cup Y_k''\) with 
\[Y_k'=\{\sigma _{k,0},\ldots,\sigma
_{k,k-1}\}=W_{\{s_1,\ldots,s_{k-1}\}}/W_{\{s_2,\ldots,s_{k-1}\}}\]
and
\[Y_k''=\{\sigma _{k,k},\ldots,\sigma
_{k,n-1}\}=w(W_{\{s_{k+2},\ldots,s_{n-1}\}}\backslash
W_{\{s_{k+1},\ldots,s_{n-1}\}}),\] where \(w=s_ks_{k-1}\cdots s_1\).

Let \(Z=W_{\{s_{k+1},\ldots,s_{n-1}\}}\),
\(Z_1=W_{\{s_{k+2},\ldots,s_{n-1}\}}\).  Then
\[
  Y_k''=w(Z_1\backslash Z),\qquad W/W_I=Y_k'\cup(Z/Z_1)w.
\]

We claim that the product
\(((Z/Z_1)w)\times Z\) is length-preserving.  Indeed, we have \((Z/Z_1)w
\subset W/W_I\) and \(Z\subset W_I\), so the claim follows from the
fact that \((W/W_I)\times W_I\) is length-preserving.

On the other hand, we have (observing \(w\) commuting with \(Z_1\)):
\[\begin{aligned}
((Z/Z_1)w)\times Z&=(Z/Z_1)\times\{w\}\times Z_1\times(Z_1\backslash
Z)\\&=(Z/Z_1)\times Z_1\times\{w\}\times(Z_1\backslash
Z)=Z\times(w(Z_1\backslash Z)).
\end{aligned}
\]
By Proposition~\ref{lengthprevseringdecompsition}, the decomposition
\(Z\times (w(Z_1\backslash Z))\) (i.e.~\(Z\times Y_k''\)) is length-preserving also.

The subset \(Y'_k\) commutes with \(Z\), and
the decompositions \(Y'_k\times Z\) and \(Z\times Y'_k\) are
length-preserving.  So \(Z\times Y_k=Z\times(Y_k'\cup Y_k'')\) is a
length-preserving decomposition and it is equal to \((Y'_k\cup
(Z/Z_1)w)\times Z= (W/W_I) \times Z=Y_n\times Z\), which is also a
length-preserving decomposition.

Applying Proposition~\ref{lengthprevseringdecompsition} to \(Z\times
Y_k=Y_n\times Z\), we get \(P_Z P_{Y_k}=P_{Y_n}P_Z\) in \(H[[u]]\).
This completes the proof.\end{proof}

\subsection{Alternating Products Over Parabolic Subgroups}

We also need a general reduction for alternating products over
parabolic subgroups, independent of type.

Given a subset $I \subseteq S$, let $W_I$ denote the parabolic
subgroup of $W$ generated by $I$. The Coxeter diagram of $(W_I, I)$ is
denoted by $\Gamma(W_I, I)$. Since \(1\in W_I\), \(P_{W_I}\) belongs to
$H[[u]]^\times$, allowing for the definition of $P_{W_I}\!^{-1}$. We
will simplify the expression
\[
\prod_{I \subseteq S} P_{W_I}\!^{(-1)^{|S \setminus I|}} \quad \text{in} \, \left(H[[u]]^\times\right)^{\text{ab}}.
\]

If $I, J \subseteq S$ are disjoint and every element of $I$ commutes with every element of $J$, then $W_{I \cup J} = W_I \times W_J$ is length-preserving. Thus, if $\Gamma(W_{J_1}, J_1)$, $\ldots$, $\Gamma(W_{J_s}, J_s)$ are the irreducible components of $\Gamma(W_I, I)$, then $P_{W_I} = P_{W_{J_1}} \cdots P_{W_{J_s}}$ in both $H[[u]]$ and $\left(H[[u]]^\times\right)^{\text{ab}}$.

A subset $I$ of $S$ is called {\it irreducible} if $\Gamma(W_I, I)$ is
connected; $I$ is called {\it dense} if each $s \in S$ is adjacent to some element in $I$ in the Coxeter diagram $\Gamma(W, S)$.

The original expression can be rewritten as
\[
\prod_{I \subseteq S} P_{W_I}\!^{(-1)^{|S \setminus I|}} = \prod_{J:
  \text{ irreducible}} (P_{W_J})^{n_J} \quad \text{in} \, \left(H[[u]]^\times\right)^{\text{ab}},
\]
where $\{n_J\}_{J:\text{ irreducible}}$ are unique integers. These exponents can be explicitly determined: 
\begin{proposition}\label{irrCoxeter}
Let $(W,S)$ be an irreducible Coxeter system, we have
\[
\prod_{I \subseteq S} P_{W_I}\!^{(-1)^{|S \setminus I|}} = \prod_{J:\text{irreducible, dense}} P_{W_J}\!^{(-1)^{|S \setminus J|}} \quad \text{in} \, \left(H[[u]]^\times\right)^{\text{ab}}.
\]
\end{proposition}

\begin{proof}
Let $J$ be an irreducible subset of $S$. Consider the set
\[
\mathcal{I} = \left\{ I \subseteq S : \Gamma(W, J) \text{ is an irreducible component of } \Gamma(W, I) \right\}.
\]
The exponent of \(P_{W_J}\) in \(\prod_{I \subseteq S} P_{W_I}\!^{(-1)^{|S \setminus I|}}\) is given by \(\sum_{I \in \mathcal{I}} (-1)^{|S \setminus I|}\).

Assume \(s_0 \in S\) is not adjacent to any element of \(J\) in \(\Gamma(W, S)\). Partition \(\mathcal{I}\) into \(\mathcal{I}_0 \sqcup \mathcal{I}_1\), where \(\mathcal{I}_0 = \{I \in \mathcal{I} : s_0 \notin I\}\) and \(\mathcal{I}_1 = \{I \in \mathcal{I} : s_0 \in I\}\). The map \(I \mapsto I \cup \{s_0\}\) gives a bijection from \(\mathcal{I}_0\) to \(\mathcal{I}_1\).
Since \((-1)^{|S \setminus I|} = -(-1)^{|S \setminus (I \cup \{s_0\})|}\), the contributions of \(P_{W_J}\) from \(\mathcal{I}_0\) and \(\mathcal{I}_1\) to \(\prod_{I \subseteq S} P_{W_I}\!^{(-1)^{|S \setminus I|}}\) cancel.

If no such \(s_0\) exists, then \(\mathcal{I} = \{J\}\) is a singleton, rendering the proposition trivially true.
\end{proof}

\section{Hyperbolic Stabilizers of Geodesic Tubes}

The construction of geodesic tubes and their stabilizers depends on
the affine root hyperplane arrangement, so we first fix the affine Weyl
group notation.

\subsection{Affine Weyl Groups}

We fix the following notation for affine Weyl
groups. Let \( \Phi \) be a reduced root system in a real vector space
\( V \). 
The corresponding Weyl group is
denoted by \( W^f \), and is identified as a group of automorphisms of
the dual space \( V^* \).

Fix a \( W^f \)-invariant inner product \( (\cdot, \cdot) \) on \( V^* \), and let \( \langle \cdot, \cdot \rangle \) denote the canonical pairing on \( V \times V^* \). Let \( \scrA \) be the affine space associated with \( V^* \), called the \textit{apartment} of \( \Phi \).

We identify an {\it affine root} \( (\alpha, k) \in \Phi^{\rm aff}:=\Phi \times \bbZ
\) with the affine function \((\alpha ,k)_\scrA\) on \( \scrA \) defined by \( (\alpha, k)_\scrA(x) := \langle \alpha, x \rangle - k \). The corresponding affine reflection is
\[ s_{\alpha, k}(x) := x - (\langle \alpha, x \rangle - k) \alpha^\vee, \]
where \( \alpha^\vee \) is the coroot of \( \alpha \).

The group generated by \( \{s_{\alpha, k} : \alpha \in \Phi, k \in \bbZ\} \) is the affine Weyl group of \( \Phi \), denoted by \( \Waff \). Its translation subgroup is \( \{t^v : v \in Q\} \simeq Q \), where \( Q \) is the \( \bbZ \)-module spanned by the coroots \( \Phi^\vee = \{\alpha^\vee : \alpha \in \Phi\} \) and \( t^v(x) := x + v \).

Choose a set of simple roots \( B \) for \( \Phi \). The
corresponding set of linear reflections \( S^f = \{s_{\alpha, 0} :
\alpha \in B\} \) satisfies the following property: \( (W^f, S^f) \)
forms a Coxeter system. A unique set \( S \) of affine reflections
exists such that \( S \supset S^f \) and  \( (W, S) \) also forms a
Coxeter system.  In fact, \((W^f,S^f)\) is then a parabolic subgroup
of \((W,S)\).  When \(S\) is used as the color set of a chamber
complex, we write it as \(S^\aff\).

For a proper subset \( I \subset S^f \), set \( B_I = \{ \alpha \in B : s_{\alpha, 0} \in I \} \). Then \( \Phi_I \) is the root system with simple roots \(B_I\). The corresponding Weyl group and affine Weyl group are denoted by \( W^f_I \) and \( \Waff_I \), respectively. Further, \( \scrA_I \) is the affine space associated with the subspace of \( V^* \) spanned by \( \alpha^\vee \) for \( \alpha^\vee \in B_I \).

The length function of \( (W, S) \) is denoted by \( \ell \).  Notice
that \( \ell|_{W^f} \) is the length function of \( (W^f, S^f) \).

For an affine root \( (\alpha, k) \), we refer to the hyperplane \(
H_{\alpha, k} = \{ v \in \scrA : \langle \alpha, v \rangle = k \} \)
as a wall. The connected components of \( \scrA\setminus
\bigcup_{(\alpha,k)\in\Phi\times\bbZ} H_{\alpha,k} \) are
{\it alcoves}.

\subsection{Geodesic Tubes} \label{geodesictube}

With the affine walls and alcoves fixed, we now isolate the tubular
regions whose stabilizers will produce the desired generators.

Inspired by \cite{KM}, we consider subsets \(\scrT\) of \(\scrA\) satisfying:
\begin{itemize}
    \item the closure of \(\scrT\) is a union of the closures of alcoves;
    \item there is an orthogonal decomposition \(\scrA = \scrL \times \scrA'\), where \(\scrL\) is an affine line, \(\scrA'\) is an affine subspace of codimension 1, and \(\scrT = \scrL \times \scrD\) for a bounded domain \(\scrD\) in \(\scrA'\).
\end{itemize}
Moreover, \(\scrT\) should be minimal with these properties.

It follows that \(\scrT\) should be bounded by walls \(H_{\alpha ,k}\)
with \(\alpha   \perp v\), where \(v\) is a non-zero tangent
vector of \(\scrL\).  For the bounded domain \(\scrD\) to be defined by half-spaces bounded by such walls, \(\{\alpha \in \Phi:
\alpha \perp v\}\) must be a root system of rank \(\dim V -1\)
(cf.~\cite[Cor.~to Prop.~VI.4]{Bou}).  By \cite[Prop.~VI.24]{Bou}, there exists a system of simple roots \(B\) of \(\Phi\), together with
\(\beta\in B\), such that \(v\) generates the line orthogonal to
\(B\setminus \{\beta\}\).  This gives the following definition.

\begin{definition}
  Given a system of simple roots \( B \)
and a root \( \beta \in B \), let \( v= \omega^\vee_\beta \) be the corresponding fundamental coweight satisfying \( \langle \alpha, v \rangle = \delta_{\alpha \beta} \) for all \( \alpha \in B \). Let \( I_\beta \subseteq S \) be the subset such that \( B_I = B \setminus \{ \beta \} \), and set
\[
\scrW_{v} = \{ H_{\alpha, k} : \alpha \in \Phi, k \in \bbZ, \langle \alpha, v \rangle = 0 \} = \{ H_{\alpha, k} : \alpha \in \Phi_{I_\beta}, k \in \bbZ \}.
\]
For any connected component \( \scrT \) of \( \scrA \setminus
\bigcup_{H \in \scrW_v} H \), we call \((\scrT,v)\) a \textit{geodesic tube of direction \( v \)}. When \(v\) is parallel to \(\omega_\beta^\vee\), the tube has \textit{type} \(\beta\). Varying
\(B\), \(\beta\), and \(\scrT\) gives all {\it geodesic tubes}.  Now
fix \(B\) once and for all.  By \cite[IV.1.10]{Bou}, for any
geodesic tube \((\scrT,v)\), there exists a unique \(\beta \in B\) such
that \(
w.(\scrT, \omega_\beta^\vee)\) is of direction \(\omega_\beta^\vee\) for some \(w \in
W\).
\end{definition}

\begin{remark}
Given a geodesic tube \((\scrT,v)\), there exists another geodesic
tube \((\scrT,-v)\). These are the only tubes of the form \(
(\scrT, v') \). Thus a geodesic tube may be viewed as a geometric shape \( \scrT \) equipped with an \textit{orientation} \( v \). When the direction is fixed, we also write \(\scrT\) for the oriented tube \((\scrT,v)\).
\end{remark}

Consider a geodesic tube \( \scrT \) with direction \(v=
\omega^\vee_\beta \). Let \( \scrA' \) be the orthogonal complement
of \( \mathbb{R}v \) in \( V^* \), regarded as an affine
subspace of \( \scrA \). Then \( \scrA' = \scrA_{I_\beta}
\). The intersections \( \{ H \cap \scrA' : H \in \scrW_v \} \)
are the walls of the affine root system \(
(\Phi_{I_\beta})^\text{aff} \). Hence \( \scrA' \cap \scrT \)
is an alcove in this affine root system, and
\( \Waff_{I_\beta} \) acts simply transitively on the geodesic
tubes of direction \( v \).

Let \( \mathfrak{y} \) denote the barycenter of \( \scrA' \cap \scrT
\). For any hyperplane \( \scrA'' \) parallel to \( \scrA' \), the intersection of \( \scrA'' \) with \( \scrL := \mathfrak{y} + \mathbb{R}\omega^\vee_\beta \) is the barycenter of \( \scrA'' \cap \scrT \). This identifies \( \scrA \) with \( \scrL \times \scrA' \), and \( \scrT \) with \( \scrL \times (\scrA' \cap \scrT) \). The line \( \scrL \) is the \textit{central line} of \( \scrT \), which will control the stabilizer of \( \scrT \).
The central line is defined similarly for geodesic tubes oriented in any direction.

\begin{center}
\begin{tikzpicture}[scale=1]
\draw 
(-1,-1)--(2,0)--(2,4)--(-1,3)--(-1,-1);
\draw (3.5,1)--(0,1);
\draw [dashed] (0,1)--(-1,1);
\draw (-1,1)--(-3,1);

\draw [yshift= 1 cm] (3.5,1)--(0,1);
\draw [yshift=1 cm, dashed] (0,1)--(-1,1);
\draw  [yshift= 1 cm] (-1,1)--(-3,1);

\draw (4.4,1.3)--(.9,1.3);
\draw [dashed] (.9,1.3)--(-1,1.3);
\draw (-1,1.3)--(-2.1,1.3);

\draw [yshift= 1 cm] (4.4,1.3)--(.9,1.3);
\draw [dashed,yshift= 1 cm] (.9,1.3)--(-1,1.3);
\draw [yshift= 1 cm](-1,1.3)--(-2.1,1.3);
\filldraw [red] (0.45,1.65) circle (1pt);

\draw [red] (3.95,1.65) -- (-2.55,1.65);
\draw  (3.5,1)--(4.4,1.3)--(4.4,2.3)--(3.5,2)--(3.5,1);
\draw [dashed, xshift=-3.5 cm] (3.5,1)--(4.4,1.3)--(4.4,2.3)--(3.5,2)--(3.5,1);
\draw [ xshift=-6.5 cm] (3.5,1)--(4.4,1.3)--(4.4,2.3)--(3.5,2)--(3.5,1);
\node at (.4,0) {{\(\scrA '\)}};
\node at (0.5,1.8)  {{$ \mathfrak{y}$}};
\node at (4.15,1.7)  {{$\scrL$}};
\node at (2.5,2.5)  {{$\scrT$}};
\end{tikzpicture}
\end{center}

\begin{remark}
An element \((\alpha,k)\) of \(\Phi_{I_\beta}\times\bbZ\) is both an affine root of
\(\Phi^{\rm aff}\) and \(\Phi_{I_\beta}^{\rm aff}\).  If the alcove
\(\scrA '\cap \scrT\) is given by \(\{(\alpha _0,k_0)_{\scrA
  '}>0,\ldots,(\alpha _m,k_m)_{\scrA '}>0\}\), then
\(\scrT\) is given by \(\{(\alpha _0,k_0)_{\scrA
  }>0,\ldots,(\alpha _m,k_m)_{\scrA }>0\}\).
\end{remark}

\begin{proposition} \label{centralline}
The central line \( \scrL \) of a geodesic tube \((\scrT,v)\) intersects at least one alcove. Moreover, if \( w \in \Waff \) fixes \( \scrL \) pointwise, then \( w \) is the identity.
\end{proposition}

\begin{proof}
It is sufficient to consider the case that  $v=\omega^\vee_\beta$.
For the first assertion, suppose \( \scrL \) does not intersect any alcove. Then \( \scrL \) is contained in the union of the walls \( \{ H_{\alpha,k} \} \). Since there are only countably many such walls and \( \scrL \) contains uncountably many points, some wall \( H_{\alpha,k} \) contains two distinct points of \( \scrL \). As \( v\) is tangent to \( \scrL \), we have \( \langle \alpha,  v \rangle = 0 \), so \( H_{\alpha,k} \in \scrW_v \). This contradicts the fact that the barycenter \( \mathfrak{y} \) is not contained in any wall in \( \scrW_v \).

The second assertion follows from the simple transitivity of the action of \( \Waff \) on alcoves.
\end{proof}

\def \Min{\mbox{Min}}

\subsection{Geometry of affine isometries}

To analyze elements stabilizing a tube, we need standard facts about
the minimal sets and displacement vectors of affine isometries.

We use the following standard facts about affine isometries of
a Euclidean space \(\scrA\).

For an affine isometry $\ga$ of $\scrA$, we define:
$$ d_\ga = \inf\{d(x,\ga x): x\in \scrA\} \quad \text{and} \quad \Min(\ga)= \{x \in \scrA: d(x,\ga x)=d_\ga\}.$$
If $\Min(\ga)$ is non-empty, $\ga$ is called {\it
  semi-simple}. In this case, if $d_\ga>0$, we call $\ga$ {\it hyperbolic}.

The following result from \cite[Ch~2, Prop~6.2 and Th~6.8]{BH} will be used.

\begin{theorem}\label{theoremBH}
Every affine isometry $\ga$ of $\scrA$ is semi-simple and $\Min(\ga)$ is an affine subspace of $\scrA$. The restriction of $\ga$ to $\Min(\ga)$ is a translation by a vector $v_{\ga}$, called the displacement vector of $\ga$. Moreover, every $\ga$-invariant affine line is contained in $\Min(\ga)$.
\end{theorem}

The displacement vector gives the following corollary.
\begin{corollary} \label{cor1}
Let $\ga$ be an affine isometry on $V$ and let $E_1$ be a $\ga$-invariant affine subspace contained in $\Min(\ga)$. Let $E_2$ be any affine orthogonal complement of $E_1$. Then $t^{-v_\ga} \ga(E_2)=E_2$.

\begin{proof} By using the unique point in \(E_1 \cap E_2\) as the
  origin, regard \(\scrA\), \(E_1\), \(E_2\) as vector spaces.
  Then $E_2$ is the usual orthogonal complement subspace of $E_1$ in \(\scrA\). By Theorem \ref{theoremBH}, $t^{-v_\ga}\ga$ acts trivially on $E_1$, hence fixes the origin. Therefore, $t^{-v_\ga}\ga$ is a linear isometry and stabilizes the orthogonal complement $E_2$ of $E_1$.
\end{proof}
\end{corollary}

\subsection{Hyperbolic stabilizers of a geodesic tube}

These displacement-vector facts now allow us to define and control the
hyperbolic stabilizer of a tube.

For an oriented geodesic tube \((\scrT,v)\), define its hyperbolic stabilizer by
\[
\Stab^h(\scrT) = \{ w \in \Waff \mid w(\scrT) = \scrT,  w \text{ is either hyperbolic or the identity} \}.
\]
Proposition~\ref{hypergen} below shows that this set is a subgroup of \( \Waff \).

\begin{proposition} \label{prop1}
Let \( w \in \Stab^h(\scrT) \) and decompose \( w = t^{v_0} w_0 \), where \( v_0 \in Q \) and \( w_0 \in W^f \).
\begin{enumerate}
\item The displacement vector \( v_w \) is parallel to \( v \).
\item The central line \( \scrL \) of \( \scrT \) is \( w \)-invariant.
\item The vectors \( v_0 - v_w \) and \( v_w \) are perpendicular.
\end{enumerate}
\end{proposition}

\begin{proof}$ $
\begin{enumerate}
\item Let \( m >0\) be such that \( w^m \) is a translation. By Theorem \ref{theoremBH}, \( w^m \) translates by \( m v_w \) and stabilizes \( \scrT \). Hence, \( v_w \) is parallel to \(v \), a tangent vector of \( \scrL \).

\item Theorem \ref{theoremBH} guarantees that \( \Min(w) \) contains a \( w \)-invariant affine line \( \scrL' \) in the direction of \( v_w \) (and thus, \( v \)). By Corollary \ref{cor1}, \( t^{-v_w}w \) stabilizes \( \scrA_I \). Since it also stabilizes \( \scrT \), it must stabilize \( \scrA_I \cap \scrT \) and fix its barycenter \( \mathfrak{y} \). This establishes that the central line \( \scrL = \mathfrak{y} + \mathbb{R}v\) is \( w \)-invariant.

\item Decompose \( v_0 = v_1 + v_2 \) such that \( v_1 \parallel v_w \) and \( v_2 \perp v_w \). As \( w \) commutes with \( t^{mv_w} \), the linear part \( w_0 \) of \( w \) commutes with \( v_1 \). Thus, \( w_0(v_1) = v_1 \) and \( w_0(v_2) \) remains perpendicular to \( v_1 \).

Noting that
\[
( w(x), v_w ) = ( w_0(x) + v_0, v_1 ) = ( w_0(x), v_w ) + ( v_1, v_w ),
\]
we find that
\[
( w^m(x), v_w ) = ( w_0^m(x) + v_0, v_1 ) = ( w_0(x), v_w ) + m ( v_1, v_w ).
\]
Setting \( x = 0 \) and using \( w^m(x) =x+ mv_w  \), we deduce \( ( m v_w, v_w ) = m ( v_1, v_w ) \), proving that \( v_1 = v_w \). Consequently, \( v_0 - v_w = v_2 \) is perpendicular to \( v_w \).
\end{enumerate}
\end{proof}

The next result gives a criterion for determining when an element of \( \Waff \) stabilizes a geodesic tube.

\begin{theorem} \label{stabilizercriterion}
For a hyperbolic element \( w \) in \( \Waff \), \( w \) stabilizes a geodesic tube of direction \( v \) if \( v_w \) is parallel to \( v \) and \( \Min(w) \) is not contained in any wall.
\end{theorem}

\begin{proof}
Decompose \( w = t^{v_0} w_0 \), where \( v_0 \in Q \) and \( w_0 \in W^f \). Let \( m \) be such that \( w^m \) translates by \( m v_w \). Consequently, \( w \) and the translation by \( v \) commute, and we obtain \( w_0 \cdot v = v \).

For any \( H_{\alpha,k} \in \scrW_v \), we have \(\langle \alpha , v \rangle = 0\) and \( w \cdot H_{\alpha,k} = H_{w_0 \cdot \alpha, k'} \) for some integer \( k' \).
Because \( \langle \alpha , v \rangle = 0 \), we deduce
\[
\langle w_0 \cdot \alpha , v \rangle = \langle \alpha, w_0^{-1} \cdot v \rangle = \langle \alpha, v \rangle = 0.
\]
Thus, \( w \cdot H_{\alpha,k} \in \scrW_v \), implying that \( w \) permutes the geodesic tubes of direction \( v \).

Given that \( \Min(w) \) is not contained in any wall by assumption, it must not be contained in the union of all walls either. Therefore, there exists some point \( y \in \Min(w) \) that lies in the interior of a geodesic tube \( (\scrT, v) \). Since \( w \cdot y = y + v_w \in \scrT \), \( w \) must stabilize \( \scrT \).
\end{proof}

\subsection{Canonical generator of $\Stab^h(\scrT)$}

Once the hyperbolic stabilizer is cyclic, its positive translation
direction singles out a canonical generator.

By Proposition \ref{prop1}, for every hyperbolic element $w$ stabilizing $(\scrT,v)$, there is a real number \( c_w \) such that \( v_w = c_w v \). Moreover, \( \mathbb{Q} v \cap Q \) is an infinite cyclic group generated by \( c_0 v \) for some rational \( c_0 > 0 \).

The translation \( t^{c_0 \omega^\vee_\beta} \) belongs to \( \Stab^h(\scrT) \). Moreover, for any \( w \in \Stab^h(\scrT) \), we have \( w^{|W|} \in \langle t^{c_0v} \rangle \) since \( w^{|W|} \) is always a translation. 

From these observations, we deduce that the set \( \{ c_w : w \in \Stab^h(\scrT) \} \) is a subset of \( \bbZ \cdot c_0 \cdot |W|^{-1} \) and is discrete. Among the elements of \( \Stab^h(\scrT) \), there exists \( w_{\min} \) such that \( c_{w_{\min}} \) is the minimal positive value.

\begin{proposition} \label{hypergen}
The set $\Stab^h(\scrT)$ is an infinite cyclic subgroup generated by $w_{\min}$.
\end{proposition}

\begin{proof}
First, we show that $\Stab^h(\scrT)$ is a subgroup. For $w, w' \in \Stab^h(\scrT)$, we show that $ww' \in \Stab^h(\scrT)$.
By Theorem \ref{theoremBH}, $ww'$ is semi-simple, and hence hyperbolic if $d_{ww'}>0$.
If $d_{ww'}=0$, we show that $ww'$ is the identity.
By Proposition \ref{prop1}, the central line $\scrL$ is both $w$-invariant and $w'$-invariant and hence is also $(ww')$-invariant.
By Theorem \ref{theoremBH}, $\scrL$ is contained in $\Min(ww')$, which implies that $\scrL$ is fixed pointwise by $ww'$.
By Proposition \ref{centralline}, $ww'$ is the identity.

For a hyperbolic element $w$ stabilizing $\scrT$, the definition of $w_{\min}$ and the division algorithm give an integer $m$ such that $c_w = m c_{w_{\min}}$. Thus, $w (w_{\min})^{-m}$ fixes $\scrL$ pointwise. By Proposition \ref{centralline}, \(w (w_{\min})^{-m}=1\), so \(w = (w_{\min})^m\). Hence $\Stab^h(\scrT)$ is cyclic and \(w_{\min}\) is a generator.
\end{proof}

We call \( w_{\min} \) the \textit{canonical generator} of the oriented geodesic tube \(\scrT\).
The canonical generator of
\(\scrT '\) is \(x w_{\min}x^{-1}\) if \(\scrT
'=x.(\scrT,v)\).
It remains to compute the displacement vector of \(w_{\min}\).

\begin{theorem} \label{displacementvector}
Let $(\scrT,v)$ be a geodesic tube of type $\beta$ with canonical generator \( w_{\min}\).
The displacement vector of \( w_{\min} \) is given by
\(
\frac{(\beta^\vee, \omega^\vee_\beta)}{(\omega^\vee_\beta, \omega^\vee_\beta)} \omega^\vee_\beta.
\)
In other words, it is the orthogonal projection of \( \beta^\vee \) onto \( \omega^\vee_\beta \).
\end{theorem}

\begin{proof}
It is sufficient to consider the case $v= \omega_\beta^\vee$.
We prove this in two steps:
\begin{enumerate}
\item For all non-identity \( w \in \Stab^h(\scrT) \), 
\(
|c_{w}| \geq \frac{(\omega^\vee_\beta, \beta^\vee)}{(\omega^\vee_\beta, \omega^\vee_\beta)}.
\)
\item There exists some \( w \in \Stab^h(\scrT) \) such that 
\(
c_{w} = \frac{(\omega^\vee_\beta, \beta^\vee)}{(\omega^\vee_\beta, \omega^\vee_\beta)}.
\)
\end{enumerate}

\textbf{Step 1:} 
Let \( w = t^{v_0} w_0 \) where \( v_0 \in Q \) and \( w_0 \in W^f \). By part (3) of Proposition \ref{prop1}, \( v_0-v_w \) and \( v_w = c_w \omega^\vee_\beta \) are orthogonal. Thus,
\begin{equation}\label{eqcw}
c_{w} = \frac{(\omega^\vee_\beta, v_{w})}{(\omega^\vee_\beta,\omega^\vee_\beta)}= \frac{(\omega^\vee_\beta, v_{w}+v_0-v_{w})}{(\omega^\vee_\beta,\omega^\vee_\beta)}= \frac{(\omega^\vee_\beta, v_0)}{(\omega^\vee_\beta,\omega^\vee_\beta)}.
\end{equation}
Express \( v_0 \) as a linear combination of simple coroots \( \sum_{\alpha \in B} n_\alpha \alpha^\vee \) with \( n_\beta \neq 0 \). Then,
\[
  |c_{w}| = \left|\frac{(\omega^\vee_\beta, v_0)}{(\omega^\vee_\beta,\omega^\vee_\beta)}\right|= \left|\frac{(\omega^\vee_\beta, n_\beta \beta^\vee)}{(\omega^\vee_\beta,\omega^\vee_\beta)} \right|\geq 
\frac{(\omega^\vee_\beta, \beta^\vee)}{(\omega^\vee_\beta,\omega^\vee_\beta)}. 
\]

\textbf{Step 2:}
Note that \( t^{-\beta^\vee}(\scrT) \) is also a geodesic tube of direction \( \omega^\vee_\beta \), and \( \Waff_{I_\beta} \) acts transitively on such tubes. Hence there is \( w' \in \Waff_{I_\beta} \) such that \( t^{-\beta^\vee}(\scrT) = w' (\scrT) \). Then \( w = t^{\beta^\vee} w' \) stabilizes \( \scrT \). Write \( w' = t^{v_0} w_0 \), with \( v_0 \in \scrA_{I_\beta}\) and \( w_0 \in W^f_{I_\beta} \). Then \( w = t^{\beta^\vee + v_0} \tilde{w_0} \), and Equation (\ref{eqcw}) gives
\[
c_w = \frac{(\omega^\vee_\beta, \beta^\vee + v_0)}{(\omega^\vee_\beta, \omega^\vee_\beta)} = \frac{(\omega^\vee_\beta, \beta^\vee)}{(\omega^\vee_\beta, \omega^\vee_\beta)}.
\]

\end{proof}

\subsection{Straight generators of geodesic tubes}

The generators needed in Theorem~\ref{main1} must also be straight, so
we now relate canonical generators of tubes to straight conjugacy
classes.

Recall that \( S^\aff \) determines a specific alcove in \( \scrA \), called the \textit{fundamental alcove}.

\begin{definition}
 A geodesic tube \( \scrT \) is \textit{fundamental} if it contains the fundamental alcove.
\end{definition}

The group \( \Waff \) acts transitively on the geodesic tubes of type \( \beta \). The set
\[
\left\{ w_{\min} : w_{\min} \text{ is the canonical
    generator of a geodesic tube of type } \beta \right\}
\]
forms a single conjugacy class in \( \Waff \). We show below that this class is straight in the sense of \cite{HN}. The minimal-length elements in this conjugacy class are straight elements. The relevant geodesic tubes are those with straight \( w_{\min} \). We recall the needed definitions and results from \cite{HN}.

\begin{definition}
An element \( w \in \Waff \) is \textit{straight} if \(
\ell(w^k) = k \ell(w) \) for all \( k \in \mathbb{N} \). A conjugacy
class of \( \Waff \) is called \textit{straight} if it contains at
least one straight element.
\end{definition}

When the canonical generator of a geodesic tube \(\scrT\) is straight,
we call it the \textit{straight generator} of \(\scrT\).

When $v$ is an element in the coroot lattice such that \(v\) lies in the closed fundamental Weyl chamber, the length of the translation $t^v$ is equal to $\langle 2 \rho, v \rangle$, where $\rho$ is the half-sum of positive roots. This gives the following criterion for straightness. For $w\in \Waff$, let $\tilde{v}_w$ be an element in the $W^f$-orbit of $v_w$ contained in the closed fundamental Weyl chamber, where $v_w$ is the displacement vector defined in Theorem \ref{theoremBH}.

\begin{proposition} \label{straightelementlength}
With this notation, $w$ is straight if and only if $w$ admits a word expression of
length $\langle 2 \rho, \tilde{v}_w \rangle$.  In this case, $w$ is of length $\langle 2 \rho, \tilde{v}_w\rangle$.
\end{proposition}

The following criterion of He and Nie \cite[Lemma 2.7]{HN} is very
useful.

\begin{theorem}  \label{theoremHN}
  Let \(w \in \Waff \) and let \(K\subset
{\rm Min}(w)\) be an affine subspace with \(w(K)=K\).  Assume that
\(K\) contains a point on an alcove \(C\).  Then
\(x^{-1}wx\) is straight, where \(x\in \Waff \) is such that
\(x\) takes the fundamental alcove to \(C\).
\end{theorem}

\begin{theorem} \label{straightlength}
A fundamental geodesic tube of type \( \beta \) with a straight generator exists for all \( \beta \in B \).  Moreover, its canonical generator has length
  \(
  \frac{(\omega_\beta^\vee, \beta^\vee)}{(\omega_\beta^\vee, \omega_\beta^\vee)} \langle 2 \rho, \omega_\beta^\vee \rangle.
  \)
\end{theorem}


\begin{proof}
Start with any geodesic tube \(\scrT\) of type \(\beta\) and let \(w\)
be the canonical generator of \(\scrT\).
By Proposition \ref{centralline} and \ref{prop1}, the central line
$\scrL$ of $\scrT$ is contained in $\Min(w)$ and $\scrL$ intersects
some alcove \(C\). By Theorem \ref{theoremHN},
\(x^{-1} w x\) is straight.  It follows that \(x^{-1}\scrT\) is
fundamental, and the canonical generator of \(x^{-1}\scrT\)
is a straight element \(w':=x^{-1}wx\).

By Theorem \ref{displacementvector}, we have $v_{w} =\frac{(\omega_\beta^\vee, \beta^\vee)}{(\omega_\beta^\vee,\omega_\beta^\vee)}\omega_\beta^\vee$, which is contained in the closed fundamental Weyl chamber. On the other hand, $v_{w'}$ and $v_{w}$ are in the same $W^f$-orbit. Thus we have $\tilde{v}_{w'} = \tilde{v}_{w} = v_{w}.$ Therefore,
$$ \ell(w') = \langle 2 \rho, \tilde{v}_{w'} \rangle= \langle 2 \rho, v_w \rangle =
\frac{(\omega_\beta^\vee, \beta^\vee)}{(\omega_\beta^\vee,\omega_\beta^\vee)}\langle 2 \rho, \omega_\beta^\vee \rangle.$$
\end{proof}

\subsection{The proof of Theorem~\ref{main1}}

The preceding theorem supplies the required straight generators; it
remains only to match their lengths with the numerical list in Table 1.

\begin{proof}  The first assertions follow from Theorem \ref{straightlength}.
The left-hand side of Equation (\ref{eqmain1})  has been computed by \cite[Section 5]{KM}.  The right-hand side is computed using Theorem \ref{straightlength}.
It remains to show that
$d_1,\ldots,d_{n}$ given in Table 1 are indeed the lengths of
$\{w_\beta\}_{\beta\in B}$ in Theorem \ref{straightlength}.
Let \(C\) be the Cartan matrix of \(\Phi\).  Note that the
value of \(\left(\frac{(\omega_\beta^\vee,
    \beta^\vee)}{(\omega_\beta^\vee,\omega_\beta^\vee)}\right)^{-1}\)
is the \(\beta\)-th diagonal entry of \(C^{-1}\), and $\langle 2 \rho,
\omega_\beta^\vee \rangle$ is the sum of the \(\beta\)-th column of
\(C^{-1}\).  These values are obtained from Table 2 in
\cite{OV} (page 295-297).  The resulting values of \(\ell(w_\beta)\) are as follows:
\[@
  A_n: &\qquad \Dynkin{A}{n}{{n+1}{n+1}{n+1}{n+1}{n+1}}\\\noalign{\medskip}
  B_n: &\qquad \Dynkin{B}{n}{{2n-1}{2n-2}{2n-3}{n+1}{2n}}\\\noalign{\medskip}
  C_n: &\qquad \Dynkin{C}{n}{{2n}{2n-1}{2n-2}{n+2}{n+1}}\\\noalign{\medskip}
  D_n: &\qquad \Dynkin{D}{n}{{2n-2}{2n-3}{2n-4}{n+2}{n+1}{2n-2}{2n-2}}\\\noalign{\medskip}
  E_6: &\qquad \Dynkin{E}{6}{{12}{11}{9}{\quad7}{9}{12}}\\\noalign{\medskip}
   E_7: &\qquad \Dynkin{E}{7}{{17}{14}{11}{\quad8}{10}{13}{18}}\\\noalign{\medskip}
  E_8: &\qquad \Dynkin{E}{8}{{23}{13}{17}{\quad9}{11}{14}{19}{29}}\\\noalign{\medskip}
  F_4: &\qquad \Dynkin{F}{4}{{8}{5}{7}{11}}\\\noalign{\medskip}
  G_2: &\qquad \Dynkin{G}{2}{{5}{3}} 
\]
This completes the
proof of Theorem~\ref{main1}.
\end{proof}

\begin{remark}  The numbers \(d_1,\ldots,d_{n}\) in Table 1 came
  from the computation in \cite[Section 5]{KM}, and are
  given in increasing order, simply because there was no other apparent way
  of indexing them.  Now we know that these numbers should be
  canonically indexed by \(B\), the set of simple roots.
\end{remark}

  Geodesic tubes of type \( \beta \) with straight generators are
  finite in number but not necessarily unique.
  Examples show that one can have more
than one fundamental geodesic tube of type \(\beta\) with straight
generators.  However,
we shall demonstrate that Conjecture \ref{conj-b} is unaffected by the particular choice of \( w_\beta \).

\begin{proposition}
  Let \(\beta\in B\).  Let \(\scrT\) and \(\scrT '\) be
  geodesic tubes of type \(\beta\) with straight generators \(w\) and \(w'\) respectively.  Then there exists \(x \in H^\times\) such
  that \(e_w=xe_{w'}x^{-1}\) in \(H\). 
\end{proposition}

\begin{proof} The key to this is He-Nie's result
  \cite[Theorem~2.9]{HN} that \(w\) and \(w'\) are strongly conjugate
  in the sense of \cite[Section~1.2]{HN} or \cite[Section~3.2.4]{GP}.  Once this
  is known, the lemma is very similar to \cite[Lemma~8.2.1]{GP}.  The proof is the same.
\end{proof}

An important implication of the straightness attributes of \( w \) and
\( w' \) is captured in the relation \( e_{w^i} = (e_w)^i =
x(e_{w'})^i x^{-1} = x e_{{w'}^i} x^{-1} \). As a direct corollary,
given \(\bbT = \{ w^i : i \geq 0 \} \) and \(\bbT' = \{ w'^i : i \geq
0 \} \),  we have \( P_{\bbT} = xP_{\bbT'}x^{-1} \), and hence \(P_{\bbT}=P_{\bbT'}\) in \((H[[u]]^\times)^{\rm ab}\). Thus the validity of
Conjecture \ref{conj-b} doesn't depend on the choice of the straight generator \( w_i \).

\subsection{Geodesic Tube as a Convex Hull}\label{S:convex}

The zeta construction below requires recognizing a tube from the
gallery generated by its straight generator; for this we record a
convex-hull characterization.

Assume in this subsection that \( \Phi \) is irreducible. In this setting, \( \scrA \) can be identified with the geometrization of the Coxeter complex associated with \( (\Waff, S) \). The convex hull within a chamber complex is understood in the sense of \cite[Example~3.133(c)]{AB}: it is the intersection of all roots containing the given chambers. To distinguish this from the conventional notion of convexity in Euclidean space, we call it the \textit{simplicial convex hull}.

We use the same notation for both a chamber complex and its geometrization and treat \( \scrA \) as a chamber complex. The chambers of \( \scrA \) correspond precisely to the alcoves \( C \) in \( \scrA \). According to our convention, we regard the closure of a geodesic tube as chamber subcomplexes of \( \scrA \).

The following characterization of convexity will be used.

\begin{proposition} \cite[Propositions~3.94 and 3.97]{AB} \label{convexAB}
  Let \( \mathcal{D} \) be a non-empty set of chambers in \( \scrA \). Then \( \mathcal{D} \) is convex if and only if it is the intersection of sets \( \{\mathcal{C}_{(\alpha,k)}\} \), where \( \mathcal{C}_{(\alpha,k)} = \{ x \in \scrA: \langle \alpha ,x \rangle \geq k \} \) is the closed half-space corresponding to the affine root \( (\alpha,k) \). Moreover, there is a unique minimal subset \( X \) of roots such that \( \mathcal{D} \) is the intersection of sets \( \{\mathcal{C}_{(\alpha,k)} : (\alpha,k) \in X\} \).
\end{proposition}

By the remark preceding Proposition~\ref{centralline}, the closure of
a geodesic tube is convex.  The next proposition characterizes a
geodesic tube as a simplicial convex hull.

\begin{proposition}\label{convex}
Let \( (\scrT, v) \) be a geodesic tube with a canonical generator \(
w \), and let \( \bar{C} \) be an arbitrary closed chamber in \( \bar{\scrT} \). For all \( d \geq 1 \), the simplicial convex hull of \( \{w^{dn}\bar{C}: n \in \bbZ\} \) in \( \scrA \) is \( \bar{\scrT} \).
\end{proposition}

\begin{proof}
  Let \(\scrD\) be the simplicial convex hull of \(\{w^{dn}\bar C: n \in \bbZ\}\) in \(\scrA\). Since \(\bar \scrT\) is convex and  \(\bar \scrT\) contains  \(\{w^{dn}\bar C: n \in \bbZ\}\), we have \(\bar \scrT \supset \scrD\). According to Proposition \ref{convexAB}, \(\scrD = \bigcap_{(\alpha, k) \in X} \scrC_{(\alpha, k)}\), where \(X\) is the unique minimal set of affine roots defining \(\scrD\).

  Without loss of generality, assume that the direction \(v\) of \(\scrT\) is \(\omega_\beta^\vee\) for some \(\beta \in B\). Let \(m > 0\) be such that \(w^{md}\) is a translation, necessarily by a positive multiple \(rv\) of \(v\). For any \((\alpha, k) \in X\), the translation of \(\scrC = \scrC_{(\alpha, k)}\) by \(rv\) is \(\scrC' = \scrC_{(\alpha, k + \langle \alpha, rv \rangle)}\). As \(\bar \scrT\) is stabilized by translation by \(rv\), \((\alpha, k + \langle \alpha, rv \rangle) \in X\) due to the uniqueness of \(X\).

  If \(\langle \alpha, rv \rangle \neq 0\), then either \(\scrC \subsetneq \scrC'\) or \(\scrC' \subsetneq \scrC\). This implies that one can remove either \(\scrC\) or \(\scrC'\) from the intersection \(\bigcap_{(\alpha, k) \in X} \scrC_{(\alpha, k)}\) and still get \(\scrD\), violating the minimality of \(X\). Hence, \(\langle \alpha, rv \rangle = 0\).

  This implies that every \(\alpha \in X\) is an affine root of the root system \(\Phi_I\), where \(I\) is such that \(B_I = B \setminus \{\beta\}\). Hence \(\bar \scrT \subset \scrD\), and therefore, \(\scrD = \bar \scrT\).
\end{proof}

\begin{proposition}  Let \((\scrT,v)\) be a geodesic tube with straight generator \(w\). 
Let \(\scrT '\) be the fundamental geodesic tube of the same
  direction \(v\).  Then \(w\) is also the canonical generator of $\scrT'$.
\end{proposition}

\begin{proof}  
  Let \(w '\) be the canonical generator of $\scrT'$.   We know that \(\scrT '=x\scrT\) for some \(x
  \in  \Waff\) such that \(x(v)=v\).  Then \(w'=xwx^{-1}\).
  Let \(d\) be such that \(w^d\) is a translation, say by \(cv\).
  Then \((w')^k\) is the translation by \(x(cv)=cv\), so
  \(w^k=(w')^k\).  Consequently, \(\scrT '\) is the simplicial
  convex hull of \(\{w^{kn}\bar{C}_{\rm fun}: n\in \bbZ\}\), where \(C_{\rm
    fun}\) is the fundamental alcove, by the preceding proposition.

The straightness of \(w\) implies that \(w C_{\rm fun}\) is on a
minimal gallery from \(C_{\rm fun}\) to \(w  C_{\rm fun}\), so
\(wC_{\rm fun}\) is in \(\scrT '\) as \(\scrT '\) is convex.  The same argument shows that
\(w^n C_{\rm fun}\) is in \(\scrT '\) for all \(n \in \bbZ\).
This implies \(w \in\Stab^h(\scrT ')\), so \(w=(w')^m\) for
some \(m > 0\).  But \(v_w=v_{w'}\) by Theorem~\ref{straightlength}, so we must have \(m=1\).
\end{proof}

\section{Zeta functions}

We formulate the zeta side through finite digraphs, then specialize
these digraphs to galleries in the building.

\subsection{Digraphs} \label{digraph}

A \textit{digraph} (or \textit{directed graph}) is a quadruple \(\scrG=(V,E,o,t)\), where \(V\) is the set of \textit{vertices}, \(E\) the set of \textit{directed edges}, and \(o,t:E\to V\) the \textit{origin} and \textit{terminus}. The digraphs considered here are \textit{locally finite}: for every \(v \in V\), both the out-degree and in-degree of \(v\) are finite:
\[
\deg_{\text{out}}(v) := |\{e \in E: o(e) = v\}| \quad \text{and} \quad \deg_{\text{in}}(v) := |\{e \in E: t(e) = v\}|.
\]

A \textit{path} on \(\scrG\) is a tuple \(p=(e_0, \ldots, e_{n-1})\) in \(E\) such that \(t(e_i) = o(e_{i+1})\) for \(i = 0, \ldots, n-2\). It has length \(n\), origin \(o(p)=o(e_0)\), and terminus \(t(p)=t(e_{n-1})\). If \(p = (e_0, \ldots, e_{n-1})\) and \(q = (f_0, \ldots, f_{m-1})\) satisfy \(t(p)=o(q)\), their \textit{concatenation} is \(pq = (e_0, \ldots, e_{n-1}, f_0, \ldots, f_{m-1})\), of length \(n+m\). Concatenation is associative: if \(pp'\) and \(p'p''\) are defined, then \(p(p'p'') = (pp')p''\).

A path \( p \) is a \textit{cycle} if \( o(p) = t(p) \). If \( c \) is a cycle of length \( n \), then \( c^m \), the concatenation of \( m \) copies of \( c \), is a cycle of length \( nm \). A cycle \( c \) is \textit{primitive} if \( c \neq (c')^m \) for any \( m > 1 \) and cycle \( c' \).

For a cycle \(c=(e_0,\ldots,e_{n-1})\), we often identify the indexing
set \(\{0,\ldots,n-1\}\) with \(\bbZ/n\bbZ\).  If
\(c'=(e_0',\ldots,e_{n-1}')\) is another cycle of the same length, we say
\(c'\) is {\it equivalent} to \(c\) if there exists \(k \in
\bbZ/n\bbZ\) such that \(c'_i=c_{i+k}\) for all \(i\).  This is
an equivalence relation.  Its equivalence classes are called
the {\it cycle classes}.  A cycle class is called {\it primitive} if one or all of
its members are primitive.

Assume that \( \scrG \) is finite, i.e., both \( V \) and \( E \) are finite sets. The \textit{zeta function} of \( \scrG \) is defined as
\[
Z_\scrG(u) := \exp\left(\sum_{m \geq 1} \frac{N_m(\scrG)}{m} u^m\right) \in \mathbb{Q}[[u]],
\]
where \( N_m(\scrG) \) is the number of cycles of length \( m \) in \( \scrG \), and \( \exp(x) = \sum_{n=0}^{\infty} \frac{x^n}{n!} \) is viewed as a formal power series.

Following \cite{MS}, we define  the {\it adjacency matrix} \(A_X\) of \(X\) is the \(V\times V\)
matrix \(A_X\) such
that 
\[
(A_X)_{xy}=|\{e \in E: o(e)=x, t(e)=y\}|.
\]
Following \cite{KS}, we define the {\it adjacency operator} \(\tilde
A_X\) on \(C(V)\), the
space of complex-valued function on \(V\), as follows:
\[
(\tilde A_X.f)(x)=\sum_{e\in E: o(e)=x} f(t(e)).
\]
If \(X\) is finite, then \(\{\delta _x\}_{x\in V}\) is a basis of
\(C(V)\) (where \(\delta _x\) is the characteristic function of
\(x\)), and the matrix of \(\tilde A_X\) with respect to this basis is
\(A_{X^{\rm op}}\), where \(X^{\rm op}:=(V,E,t,o)\) is the {\it
  opposite digraph}.  The matrix \(A_{X^{\rm op}}\) is equal to the
transpose of \(A_X\).

\begin{theorem}\label{zetadg}
  Let \(X\) be a finite digraph.  Then
\[
Z_X(u) = Z_{X^{\rm op}}(u)=\det(I-A_Xu)^{-1}=\det(I-\tilde A_X u)^{-1}=\prod_{[C]} (1-u^{l([C])})^{-1},
\]
where the product is over all primitive cycles classes of \(X\), and
\(l([C])\) is the length of any member of \([C]\).
\end{theorem}

This is Lemma 2.2 and Theorem 2.3 of \cite{KS} (see also \cite[Theorem
4]{MS}), where the
authors assumed that \(X\) is strongly connected.  The same argument
shows that the theorem is valid without this assumption.

Let \( \Gamma \) act on \( \scrG \) by digraph automorphisms, and denote the \( \Gamma \)-orbit of \( x \) by \( \ringO(x) \). The formulas \( \bar{o}(\ringO(e)) = \ringO(o(e)) \) and \( \bar{t}(\ringO(e)) = \ringO(t(e)) \) define maps \( \bar{o}, \bar{t} : \Gamma \backslash E \to \Gamma \backslash V \). Thus \( (\Gamma \backslash V, \Gamma \backslash E, \bar{o}, \bar{t}) \) is a digraph, denoted \( \Gamma \backslash \scrG \), the \textit{quotient} of \( \scrG \) by \( \Gamma \).

\subsection{Digraph of \(\bfw\)-galleries}\label{bfw-gallery}

The digraphs relevant to our zeta functions are obtained by taking
chambers as vertices and galleries of fixed affine Weyl type as
directed edges.

Let \(G\) be a split, simple, simply connected algebraic group over a
non-archimedean local field \(F\) and let \(Z\) be its center.  Let
\(\scrB G\) be the Bruhat-Tits building of \(G\).  Then
\(\scrB G\) is a colored chamber complex, colored by \(S^\aff\), which
is associated to the root system \(\Phi\)  of \(G\), \(G(F)\) acts
on \(\scrB G\) by color-preserving simplicial automorphisms, and
\(Z(F)\) is the kernel of the action.
Here \(S^\aff\) denotes the affine simple reflection set \(S\), viewed
as the set of colors of the building.

Let \(\tilde \Gamma\) be a discrete subgroup of \(G(F)\) such that its
image \(\Gamma\) in \(G(F)/Z(F)\) acts on the set of vertices of
\(\scrB G\) freely.  It follows that \(\Gamma\) acts on \(|\scrB G|\),
the geometrization of \(\scrB G\), freely and \(|\scrB G|\to
\Gamma\backslash |\scrB G|\) is a covering map.  We assume further
that \(\Gamma\backslash|\scrB G|\) has a (unique) simplicial structure
such that \(|\scrB G|\to
\Gamma\backslash |\scrB G|\) is the geometrization of a simplicial
map.  This assumption amounts to the following: let
\(\{v_0,\ldots,v_k\}\), \(\{w_0,\ldots,w_k\}\) be two simplices of
\(\scrB G\) such that \(v_i\) and \(w_i\) are in the same
\(\Gamma\)-orbit for \(i=0,\ldots,k\), then there exists \(\gamma \in
\Gamma\) such that \(\gamma.v_i=w_i\) for all \(i\).  For example, we
can have \(\Gamma=\{1\}\).

We denote by \(\Gamma\backslash \scrB G\) the simplicial complex whose
existence is the assumption we just made.  Then
\(X:=\Gamma\backslash\scrB G\) is a colored chamber complex.
Therefore, a gallery on \(X\) has a {\it type} (\cite[3.22]{AB}) which
is a finite sequence in \(S^\aff\).  We identify such a sequence
as a {\it word} in \(S^\aff\).

\begin{definition} Let \(w\) be a non-trivial element of \(\Waff\),
  and let \(\bfw\) be a reduced word in \(S^\aff\) representing
  \(w\).  We define a digraph \(\scrG_\bfw(X)\) as follows: the set of
  vertices is \({\rm Ch}(X)\),  the set of chambers of \(X\), and the directed edges are
  galleries on \(X\) of type \(\bfw\) (to be called \(\bfw\)-galleries
  for short).  If \(e=(C_0,\ldots,C_d)\) is a \(\bfw\)-gallery,
  \(o(e)=C_0\) and \(t(e)=C_d\).  We call \(\scrG_\bfw(X)\) the {\it
    digraph of \(\bfw\)-galleries} of \(X\).
\end{definition}

Here \(\delta\) denotes the Weyl distance on the building.

\begin{lemma} \label{w-BG} Let \(C,C'\) be two chambers on \(\scrB G\).  There is a
  \(\bfw\)-gallery from \(C\) to \(C'\) exactly when
  \(\delta(C,C')=w\).  If \(\delta(C,C')=w\), the \(\bfw\)-gallery
  from \(C\) to \(C'\) is unique.  Consequently, \(\scrG_\bfw(\scrB
  G)\) has no multiple directed edges, and the in-degree and
  out-degree of every vertex is \(q^{\ell(w)}\).
\end{lemma}

\begin{proof} The first two statements follow from \cite[Propositions~4.81,
  4.83]{AB}.  By \cite[6.17]{AB}, the in/out-degrees of every vertex
  is \(|IwI/I|\), where \(I\) is a suitable Iwahori subgroup.  This
  number is computed in \cite[3.3.1]{Tits}, and in our case it
  is simply \(q^{\ell(w)}\).
\end{proof}

\subsection{Digraph of \(w\)-galleries}\label{w-gallery}

The digraph \(\scrG_{\bfw}(X)\) depends only on \(w\), rather than on the reduced word \(\bfw\). Accordingly, we also denote
\(\scrG_{\bfw}(X)\) by \(\scrG_w(X)\).

\begin{lemma}\label{lift-gallery}
  Let \(e\) be a \(\bfw\)-gallery on \(X\), from \(C\) to
  \(C'\), and let \(\tilde C\) be a chamber of \(\scrB G\) lying above
  \(C\).  There is a unique \(\bfw\)-gallery \(\tilde e\) on \(\scrB
  G\) starting from \(\tilde C\) lying
  above \(e\).
\end{lemma}

\begin{proof} By standard covering space theory, there is a unique
  gallery \(\tilde e\) starting from \(\tilde C\) and lying above
  \(e\).  This \(\tilde e\) is a \(\bfw\)-gallery because \(\scrB G\to
  X\) is color-preserving.
\end{proof}

\begin{lemma} \label{freeq-wdg}
  The group \(\Gamma\) acts on the vertices of
  \(\scrG_\bfw(\scrB G)\) freely.  The quotient digraph
  \(\Gamma\backslash \scrG_\bfw(\scrB G)\) can be identified with \(\scrG_\bfw(X)\).
\end{lemma}

\begin{proof}  The group \(\Gamma\) acts freely on the chambers of
  \(\scrB G\).  There is a natural map from
  \(\pi:E(\scrG_\bfw(\scrB G))\to E(\scrG_\bfw (X))\).  It is
  surjective according to the preceding lemma.  It remains to
  show that each fiber of \(\pi\) is a \(\Gamma\)-orbit.  Let
  \(e,e'\in E(\scrG_\bfw (\scrB G))\) be such that
  \(\pi(e)=\pi(e')\).  Write \(e=(C_0,\ldots,C_d)\),
  \(e'=(C_0',\ldots,C_d')\).  Then \(C_0'=\gamma.C_0\) for some
  \(\gamma\in \Gamma\).  By the preceding lemma again, we have \(e'=\gamma.e\).
\end{proof}

\begin{corollary}  The digraph \(\scrG_\bfw(X)\) depends only on
  \(w\), not on the particular reduced word \(\bfw\) representing
  \(w\).  The digraph \(\scrG_\bfw(X)\) has in-degree and out-degree
  \(q^{\ell(w)}\) at every vertex.
\end{corollary}

\begin{proof} The statements are true in the special case
  \(\Gamma=\{1\}\) by Lemma~\ref{w-BG}.  The general case then follows
  from this special case by Lemma \ref{freeq-wdg}.
\end{proof}

\subsection{Adjacency Operator of Digraphs}\label{adjacency matrix w-gallery}

We next identify the adjacency operator of \(\scrG_w(X)\) with the
Hecke operator \(e_w\), which converts the digraph zeta function into
a determinant of a Hecke operator.

Recall that \( G(F) \) acts transitively on the chambers of \( \scrB G
\). Fix a chamber \( C_{\text{fun}} \) in \( \scrB G \) and define \(
I = \Stab_{G(F)}(C_{\text{fun}}) \), the corresponding Iwahori
subgroup. Recall also that the adjacency operator \( \tilde{A}_{\scrG_w(X)} \)acts on the space
\[
\bbC({\rm Ch}(X)) =\bbC(\Gamma \backslash {\rm Ch}(\scrB G)) = \bbC(\Gamma \backslash G(F)/I),
\]
 the space of \( \mathbb{C} \)-valued locally
constant functions on \( G(F) \) that are left \(\Gamma\)-invariant and right \( I \)-invariant. 

Now suppose \( C \) lifts to the vertex \( x C_{\text{fun}} \) in \(
\scrG_w(X) \). We have \( (C, C') \in E \) if and only if there is a
lifting \(yC_{\rm fun}\) of \(C'\) in \(\scrG_w(X)\) such that \(
x^{-1} y \in IwI \). Let \( IwI = \coprod_{i=1}^{N} g_i I \), where \(
N = q^{l(w)} \).
Thus there are \( N \) directed edges out of \( \Gamma x C_{\text{fun}} \), with termini \( \Gamma x g_i C_{\text{fun}} \). The action of \( \tilde{A}_{\scrG_w(X)}  \) is
\[
(\tilde{A}_{\scrG_w(X)} .f)(x) = \sum_{i=1}^{N} f(x g_i).
\]

The Iwahori-Hecke algebra \( H(G(F), I) \) identifies with the Hecke algebra \( H \) of \( (\Waff, S^\aff) \) over \( \mathbb{C} \), with parameter \( q \), the cardinality of the residue field of \( \mathcal{O}_F \). Under this identification, the characteristic function of \( IwI \) corresponds to \( e_w \). The vector space \(\bbC({\rm Ch}(X))  \) is an \( H(G(F), I) \)-module, where \( \varphi \in H(G(F), I) \) acts on \( f \in \bbC({\rm Ch}(X))  \) by
\[
(\varphi.f)(x) = \int_G \varphi(g) f(xg) \, dg,
\]
where \( dg \) is the Haar measure on \( G(F) \) such that \( I \) has measure \( 1 \).
Then,
\[
(e_w.f)(x) = \int_G e_w(g) f(xg) \, dg = \sum_{i=1}^{N} f(x g_i) = (\tilde{A}_{\scrG_w(X)} .f)(x).
\]

Therefore, Theorem \ref{zetadg} gives the following:

\begin{theorem}\label{zetadg2}
Suppose that \( X \) is finite. The zeta function of \(\scrG_w(X)\) satisfies
\[
Z_{\scrG_w(X)}(u) = \det(1 - e_w u |  \bbC({\rm Ch}(X)) )^{-1}.
\]
\end{theorem}

\subsection{The case of a straight \(w\)} \label{straight-case}

When \(w\) is straight, powers of \(e_w\) match the Hecke basis
elements \(e_{w^n}\), so the zeta determinant becomes a Poincar\'e
series for the monoid generated by \(w\).

Now assume that \(w\) is straight.  By a {\it path of \(\bfw\)-galleries on
\(X\)}, we mean a path on \(\scrG_\bfw(X)\).  Therefore, it is a sequence
\(e_0,\ldots,e_n\) of \(\bfw\)-galleries.  Furthermore, if we write
\(e_i=(C_{i0},\ldots,C_{id})\), then \(C_{id}=C_{i+1,0}\).  It is
clear that
the gallery
\((C_{00},\ldots,C_{0d},C_{11},\ldots,C_{1d},\ldots,C_{n1},\ldots,C_{nd})\)
is a \(\bfw^n\)-gallery, which we term {\it
the associated \(\bfw^n\)-gallery}.  Notice that \(\bfw^n\) is
reduced by our assumption.

\begin{proposition} \label{str-apt}
  A path of \(\bfw\)-galleries of length \(n\) on \(\scrB G\) lies in an apartment.
  More precisely, if \((C_0,\ldots,C_N)\) is the \(\bfw^n\)-gallery
  associated to that path, then there exists an apartment of \(\scrB G\)
  containing \(C_0,\ldots,C_N\).
  A path of \(\bfw\)-galleries from \(C\) to \(C'\) exists if and only if
  \(\delta(C,C')\) is a positive power of \(w\).  When such a path
  exists, it is unique.  Consequently, \(\scrG_w(\scrB G)\) is a
  directed forest.
\end{proposition}

\begin{proof} The first statement follows from \cite[5.77]{AB}.  The
  second statement follows from Lemma~\ref{w-BG}.  The final statement
  says that \(\scrG_w(\scrB G)\) has no non-trivial cycles, and it is
  clear from the fact that \(\bfw ^n\) is a reduced word
  representing \(w^n\), which is necessarily non-trivial for \(n \geq 1\).
\end{proof}

\begin{remark}\label{inf-path}
  Let \(Y\) be a digraph, and \(a,b \in
  \{-\infty\}\cup\bbZ\cup\{\infty\}\) with \(a < b\).  A sequence
  \(\{e_i : a < i < b\}\) is called a {\it path} on \(Y\) if
  \((e_s,\ldots,e_t)\) is a path on \(Y\) for all integers \(a < s
  \leq t< b\).  When \(a=-\infty\) or \(b=\infty\) (or both), this
  extends the definition given in section~\ref{digraph} by allowing
  infinite paths.  The first statement of the preceding proposition is valid for
  an infinite path by the same proof.
\end{remark}

\begin{proposition} \label{zeta-P_bbT}
Suppose that \(X\) is finite.  The zeta function of
  \(\scrG_w(X)\) satisfies
\[
Z_{\scrG_w(X)}(u^{\ell(w)})=\det(P_{\bbT} \mid \bbC({\rm Ch}(X))),
\]
where \(\bbT\) is the subset \(\{w^ i : i \geq 0\}\) of \(\Waff\).
\end{proposition}

\begin{proof} By Theorem~\ref{zetadg2},
\(Z_{\scrG_w(X)}(u^{\ell(w)})\) is the determinant of the operator
  \((1-e_w u^{\ell(w)})^{-1}\) on \(\bbC({\rm Ch}(X))\).  But
  \[
    (1-e_w u^{\ell(w)})^{-1}=\sum_{n\geq0} e_w\!^n u^{n\ell(w)} \in \End(\bbC({\rm Ch}(X)))[[u]],
  \]
  and we have \(e_w\!^n=e_{w^n}\) by the straightness of \(w\).  Thus
  the right-hand side is exactly the Poincar\'e series \(P_{\bbT}\).
\end{proof}

\subsection{Circular Geodesic Tubes Zeta function}

We now reinterpret cycles in the gallery digraph as closed, or
circular, geodesic tubes in the quotient complex.

Now assume that \(w\) is straight and it is the canonical generator of
a geodesic tube \((\scrT,v)\) on \(\scrA\) of type \(\beta\) such that
\(\scrT\) contains the fundamental alcove \(C_{\rm fun}\). We also
fix an identification of \( \scrA \) with a specific apartment of \(
\scrB G \). Put \(d=\ell(w)\). As in Section~\ref{S:convex}, the same notation is often used for both a chamber complex and its geometrization.

\begin{definition} \label{def-circ}
  A \textit{circular geodesic tube of type \( \beta \)} on \( X \) is a pair \( (c, N) \), where \( N \geq 1 \) and \( c: \bar{\scrT} \to X \) is a simplicial map satisfying:
  \begin{enumerate}
    \item[(i)] On each chamber, \(c\) restricts to a color-preserving simplicial isomorphism onto a chamber of \(X\),
    \item[(ii)] The geometric realization of \(c\) is locally injective,
    \item[(iii)] \( c \) can be factored through \( \bar{\scrT} \to \langle w^N \rangle \backslash \bar{\scrT} \).
  \end{enumerate}
  The integer \(N\) is a specified period, not necessarily the least positive period. The integer \( \ell(w) N \) is the \textit{length} of \( (c, N) \).
\end{definition}

The \textit{circular geodesic tube zeta function} of \( X \) of type \( \beta \) is defined as
\[
Z^{(\beta)}_X(u) := \exp\left(\ell(w)\sum_{m \geq 1} \frac{N_m^{(\beta)}(X)}{m} u^m\right) \in \mathbb{Q}[[u]],
\]
where \( N_m^{(\beta)}(X) \) is the number of circular geodesic tubes of type \( \beta \) of length \( m \) in \( X \). These are counted as maps with specified period, before identifying maps that differ by precomposition with a power of \(w\). In particular, \(N_m^{(\beta)}(X)=0\) unless \(d\mid m\), and equivalently
\[
Z^{(\beta)}_X(u)=\exp\left(\sum_{N\geq1}
\frac{N_{dN}^{(\beta)}(X)}{N}u^{dN}\right).
\]

\subsection{The proof of Theorem \ref{main3}}

We prove Theorem~\ref{main3}, with \(Z^{(\beta)}_X(u)\) as defined above.

It remains to identify circular geodesic tubes of length \( \ell(w) N \) with cycles of \( \bfw \)-galleries of length \( N \). This gives \(N_{dN}^{(\beta)}(X)=N_N(\scrG_w(X))\), so the normalization above and Proposition \ref{zeta-P_bbT} give
\[
Z^{(\beta)}_X(u) = Z_{\scrG_w(X)}(u^{\ell(w)})=\det(P_{\bbT} \mid
\bbC({\rm Ch}(X)))\qquad \mbox{, where }\bbT=\{w^N:N\geq 0\}.
\]

\( \bfw \)-galleries on the Coxeter complex \( \scrA \) are defined
as in Section~\ref{w-gallery}. By \cite[Exercise~3.7 after 3.86]{AB}, there is a unique infinite (cf.~\ref{inf-path}) path \( \{e_i\}_{i\in\bbZ} \) of \( \bfw \)-galleries on \( \scrA \), with \( e_i \) going from \( w^i.C_{\rm fun} \) to \( w^{i+1}.C_{\rm fun} \), and this path lies in \( \scrT \). If \( (c,N) \) is a circular geodesic tube of type \( \beta \) on \( X \), then \( \{c_*(e_i)\}_{i\in \bbZ/N\bbZ} \) is a cycle of \( \bfw \)-galleries on \( X \). This gives a map
\[
\phi: \left\{ \text{circular geodesic tubes of type } \beta \text{ and
    length } \ell(w)N \text{ on } X \right\}\\
\qquad \to \left\{ \text{cycles of } \bfw \text{-galleries of graph length } N \text{ on } X \right\}.
\]

To prove Theorem \ref{main3}, it remains to show the following.

\begin{theorem}
The map \( \phi \) is a bijection.
\end{theorem}

\begin{proof}
  We first record a consequence of Proposition~\ref{convex}. If two
  locally injective, color-preserving chamber maps from
  \(\bar\scrT\) to \(\scrB G\) agree on every
  \(w^iC_{\rm fun}\), they agree everywhere. Indeed, local
  injectivity prevents folding across a panel, so minimal galleries
  in \(\bar\scrT\) map to galleries of the same reduced type.
  Uniqueness of a gallery with fixed endpoints and reduced type
  makes the chamberwise agreement set simplicially convex. It
  contains the orbit \(\{w^iC_{\rm fun}:i\in\bbZ\}\), whose
  simplicial convex hull is \(\bar\scrT\).

  To construct an inverse for \( \phi \), start with a cycle of \( \bfw \)-galleries \( \{ \lambda_i \}_{i \in \bbZ/N\bbZ} \) on \( X \). Concatenating \( \lambda_0, \ldots, \lambda_{N-1} \) gives a gallery \( C_0, \ldots, C_{Nd} \) of length \( Nd \). By the discussion preceding Proposition~\ref{str-apt}, this is a \( \bfw^N \)-gallery. Extend it periodically to \( \{ C_i \}_{i \in \bbZ} \), with \( C_i \) depending only on \( i \bmod Nd \).

  Let \(\hat C_0\) be a chamber of \( \scrB G \) that lifts \( C_0
  \). As in Lemma~\ref{lift-gallery}, this periodically extended gallery has a unique lift \( \{ \tilde{C}_i \}_{i \in \bbZ} \)
  starting at \(\hat C_0\).

By Proposition~\ref{str-apt} and the remark following it, there is an
apartment \(\scrA '\) of \(\scrB G\) such that \(\scrA '\) contains
\(\tilde C_i\) for all \(i \in \bbZ\).  Choose \(g \in G(F)\)
such that \(\scrA '=g.\scrA\) and \(\tilde C_0=g.C_{\rm fun}\).

  Define \(c:\bar\scrT\to X\) by \(c(x)=\pi(g.x)\), where
  \(\pi:\scrB G\to X\) is the natural covering map. Condition (i)
  holds, and condition (ii) follows because \(g|_{\bar\scrT}\) is
  an embedding and \(\pi\) is locally injective. To verify (iii),
  let \(\gamma\in\Gamma\) carry \(\tilde C_0\) to
  \(\tilde C_{Nd}\). Uniqueness of gallery lifting gives
  \(\gamma\tilde C_j=\tilde C_{j+Nd}\) for all \(j\in\bbZ\).
  Hence \(\gamma g|_{\bar\scrT}\) and
  \(gw^N|_{\bar\scrT}\) agree on every \(w^iC_{\rm fun}\).
  The preceding observation extends this equality to the whole tube.
  Applying \(\pi\) gives \(c\circ w^N=c\), proving (iii).

The construction of \(c\) apparently depends on two choices: the
  choice of \(\hat C_0\) and of \(g\). However, if we choose \(g'\)
  in place of \(g\), then \(g\) and \(g'\) agree on every chamber
  \(w^iC_{\rm fun}\), since both carry the model gallery to the
  same lifted gallery. The preceding observation shows that their
  restrictions to \(\bar\scrT\) agree, so the resulting map \(c\)
  is unchanged. If we choose \(\hat C_0'\) in place of \(\hat C_0\),
then \(\hat C_0'=\gamma.\hat C_0\) for some \(\gamma \in \Gamma\), 
and consequently \(\tilde C_i'=\gamma.\tilde C_i\) for all \(i\) and we
may and do assume \(g'=\gamma g\), and hence the resulting map \(c\)
is unchanged.  

Thus \((c,N)\) depends only on \(\{\lambda_i\}_{i\in
  \bbZ/N\bbZ}\), not on the choice of \(\hat C_0\) or \(g\).  We
write \((c,N)=\psi(\{\lambda_i\}_{i\in \bbZ/N\bbZ})\).  It remains to
show that \(\phi\) and \(\psi\) are inverse to each other.

When \((c,N)=\psi(\{\lambda_i\}_{i\in \bbZ/N\bbZ})\),
\(\{g.e_i\}_{i\in\bbZ}\) is an infinite path of \(\bfw\)-galleries on
\(\scrA '\).  By \cite[Exercise~3.7 after 3.86]{AB} again, this path is
unique, so its associated gallery is identical to \(\{\tilde
C_i\}_{i\in\bbZ}\).  This shows that \(\{c_*(e_i)\}_{i\in\bbZ/N\bbZ}\)
is identical to \(\{\lambda_i\}_{i\in\bbZ/N\bbZ}\).  In other words,
\(\phi\circ\psi\) is the identity map, and hence \(\phi\) is surjective.  The theorem now
follows from the following lemma.
\end{proof}

\begin{lemma} The map \(\phi\) is injective.
\end{lemma}

\begin{proof}  Suppose \((c,N)\) and \((c',N)\) are such that
  \(\phi(c,N)=\phi(c',N)\).
  Then
\(c(w^i.C_{\rm fun})=c'(w^i.C_{\rm fun})\) for all \(i \in \bbZ\). 
Let \(\check C_0\) be a chamber of \(\scrB G\) lying above \(c(C_{\rm
fun})=c'(C_{\rm fun})\).  Then since \(\bar \scrT\) is simply
connected, \(c\) lifts to a unique map \(\tilde c: \bar \scrT\to \scrB
G\) such that \(\tilde c(C_{\rm fun})=\check C_0\).  Similarly \(
c'\) lifts to a unique map \(\tilde c': \bar \scrT\to \scrB
G\) such that \(\tilde c'(C_{\rm fun})=\check C_0\).  The maps \(\tilde{c}\) and \(\tilde c'\) are both locally injective and color-preserving, and both
(the concatenation of)
\(\{\tilde c(e_i)\}_{i\in\bbZ}\) and \(\{\tilde c'(e_i)\}_{i\in\bbZ}\)
lift the same (infinite) gallery on \(X\) while lifting \(c(C_{\rm
  fun})=c'(C_{\rm fun})\) to \(\check C_0\), so they are identical.
In particular, we have \(\tilde c(w^i.C_{\rm fun})=\tilde
c'(w^i.C_{\rm fun})\) for all \(i \in \bbZ\).

Let \(\scrD=\{C \in {\rm Ch}(\bar\scrT): \tilde c(C)=\tilde
c'(C)\}\).  We claim that this set of chambers is (simplicially)
convex in \(\bar\scrT\).  In other words (see \cite[3.92]{AB}), assuming that
\(\scrP=(C_0,\ldots,C_m)\) is a minimal gallery in \(\bar\scrT\) with
\(C_0,C_m\) in \(\scrD\), we claim that \(C_i\) is in \(\scrD\)
for \(i=0,\ldots,m\).  Indeed, since \(\bar \scrT\) is simplically
convex in \(\scrA\), \(C_0,\ldots,C_m\) is also a minimal gallery in
\(\scrA\), and hence its type is a reduced word in \(S^\aff\) by \cite[4.41]{AB}.  
As \(\tilde{c}\) and \(\tilde c'\) are locally injective and color-preserving,
the types of \(\tilde{c}(\scrP)\) and \(\tilde c'(\scrP)\) are the same reduced word as \(\scrP\), 
 \(\tilde{c}(\scrP)\) and \(\tilde c'(\scrP)\) are minimal galleries in \(\scrB G\)
by \cite[4.41]{AB} again.  But such a minimal gallery in \(\scrB G\) from
\(\tilde c(C_0)=\tilde c'(C_0)\) to \(\tilde c(C_m)=\tilde c'(C_m)\)
is unique by \cite[4.42]{AB}, we conclude \(\tilde c(C_i)=\tilde
c'(C_i)\) for all \(i\).  That is: \(C_i\) is in \(\scrD\) for all \(i\).
We have completed the proof that \(D\) is simplicially convex.

Moreover, \(\scrD\) contains \(\{ w^i.C_{\text{fun}} : i \in \bbZ \}\), whose simplicial convex hull is \(\scrT\) by Proposition~\ref{convex}. Hence \(\scrD\) contains \({\rm Ch}(\bar{\scrT})\). Thus \(\scrD = {\rm Ch}(\bar{\scrT})\), \(\tilde{c} = \tilde{c'}\), and \(c = c'\).
\end{proof}

\section{Results for Type $\tilde{A}_{n-1}$}

The remaining two sections prove
Conjecture~\ref{conj-b} for the affine Weyl groups of types
\(\tilde{A}_{n-1}\) and \(\tilde{C}_{n}\).  Recall that the
conjecture is the following formula:
\[
\prod_{I \subset S} P_{W_I}\!^{(-1)^{|S\setminus I|}} = \prod_{\beta \in B} P_{\bbT_\beta} \qquad \text{in} \qquad \left(H[[u]]^\times\right)^{\text{ab}}.
\]
We begin with type \(\tilde{A}_{n-1}\).

The proof proceeds by decomposing the affine Weyl group into a product
of subsets.  The construction is carried out through certain regions
of the apartment that are unions of closed alcoves.

\subsection{Setting}\label{AnSetting}

The proof in type \(\tilde A_{n-1}\) is explicit, so we first choose
coordinates in which the straight generators and fundamental regions
can be written down.

To study the case $\tilde{A}_{n-1}$ (with \(n\geq 2\)), let $\scrA$
denote the affine subspace of $\mathbb{R}^n$ consisting of vectors
$\vec{\xi} = (\xi_0, \xi_1, \ldots, \xi_{n-1})$ satisfying
$\sum_{i=0}^{n-1} \xi_i = 0$. For each $i = 0, \ldots, n-1$, define an
affine function \(x_i\) on \(\scrA\) by \(x_i(\vec{\xi}) = \xi_{i}\).

Then, there is a root system \(\Phi\) of type \(A_{n-1}\) whose set of
roots is  $\{x_i-x_j : 0 \leq i \neq j \leq n-1\}$.  The set
\[
\{a_1 := x_1 - x_0, \ldots, a_{n-1} := x_{n-1} - x_{n-2},a_n:=1+x_0-x_{n-1}\}
\]
comprises affine functions that form a system of simple affine roots.  The
corresponding closed fundamental alcove is specified by
\[
\{\vec\xi \in \scrA : x_0(\vec\xi) \leq x_1(\vec\xi) \leq \cdots \leq x_{n-1}(\vec\xi) \leq x_n(\vec\xi)\},
\]
where we adopt the convention \( x_{i+n} = x_i + 1 \) for all \( i \in
\mathbb{Z} \).  We also represent the closed fundamental alcove as
\[
 \{ x_0 \leq x_1 \leq \cdots \leq x_{n-1} \leq x_n \}
\]
A similar notation will be used for regions defined by inequalities.

Moreover, the fundamental coweights are given by
\[
\omega_i^\vee = (\underbrace{\frac{-(n-i)}{n}, \ldots, \frac{-(n-i)}{n}}_{i\text{-times}}, \underbrace{\frac{i}{n}, \ldots, \frac{i}{n}}_{(n-i)\text{-times}})
\]
for $i = 1,\ldots,n-1$. The sum of positive roots is
\[
2 \rho = (-n+1)x_0 + (-n+3)x_1 + \cdots + (n-3)x_{n-2} + (n-1)x_{n-1}.
\]

The closed fundamental alcove has vertices $v_1, \cdots, v_n$,
where \(v_n\) is the origin, and \(v_i=\omega_i^\vee\) for
\(i=1,\ldots,n-1\).  We denote
the simple reflection corresponding to \(a_i\) by \(s_i\) and put
\(S=\{s_1,\ldots,s_n\}\).  It is useful to
introduce the following convention: \(v_i:=v_{i\bmod n}\), \(s_i:=s_{i \bmod n}\), for all \(i
\in \bbZ\), where \(i\bmod n\) is the unique integer \(j\) such that
\(i\equiv j\pmod{n}\) and \(1\leq j\leq n\).  

For any \(i\in\bbZ\), the stabilizer
of $v_i$ in  the affine Weyl group \(W\), denoted by \(\Waff_{v_i}\), is the parabolic
subgroup of \((W,S)\) generated by
\(S\setminus \{s_i\}\), and is isomorphic to
\(S_n\).  It is useful to make the isomorphism explicit as follows:
\(\Waff_{v_i}\), being a subgroup of \(W\), acts on the set of affine
functions.  The subset \(\{x_i,\ldots,x_{i+n-1}\}\) is
\(\Waff_{v_i}\)-stable and the resulting homomorphism \(\Waff_{v_i}\to
S_{\{x_i,\ldots,x_{i+n-1}\}}\) is an isomorphism.  This is a useful way
  to specify certain elements of \(W\) (when they lie in one of
  the \(\Waff_{v_i}\)'s).  For example, \(s_j\) is \((x_{j-1},x_j)\) in
  \(\Waff_{v_i}\) for any \(j=i+1,\ldots,i+n-1\).

\subsection{Geodesic tube with straight generator}

We first construct the straight generators \(w_i\) supplied by
Theorem~\ref{main1} in these coordinates.

For \(1\leq i\leq n-1\), let 
\[
w_i:=(s_ns_{n-1}\cdots s_{i+1})(s_1s_2\cdots s_i).
\]
Equivalently, this is the composition of \((x_n x_{n-1}\ldots x_i)\) in \(\Waff_{v_1}\) with \((x_0x_1\ldots x_i)\) in \(\Waff_{v_0}\), and it sends \((\xi_0,\ldots,\xi_{n-1})\) to \((\xi_{i-1}-1,\xi_0,\ldots,\xi_{i-2},\xi_{i+1},\ldots,\xi_{n-1},\xi_i+1)\).
Then \(w_i^{i(n-i)}\) is translation by \(n\omega_i^\vee\) and
\(v_{w_i}=\frac{n}{i(n-i)}\omega_i^\vee\).  In particular, \(w_i\) is hyperbolic.

The set
\[
\scrT_i:=\{x_0 < x_1 < \cdots < x_{i-1} < x_0+1, x_i <
x_{i+1} <\cdots<x_{n-1}< x_i+1\}
\]
is a geodesic tube of direction \(\omega_i^\vee\), and direct verification shows that \(\scrT_i\) is \(w_i\)-invariant.  Hence \(w_i\in \Stab^h(\scrT_i)\). Here \(a_i^\vee\) denotes the coroot corresponding to the finite simple root underlying the affine simple root \(a_i\). From the
definition formula, we have \(\ell(w_i)\leq n\).  By
Theorem~\ref{straightlength}, \(\ell(w_i)\geq
\frac{(\omega_i^\vee,a_i^\vee)}{(\omega_i^\vee,\omega_i^\vee)}\langle
2\rho,\omega_i^\vee\rangle =n\).  Thus we conclude that \(w_i\) is the
canonical generator of \(\scrT_i\), and it is straight.

\subsection{The plan of the proof}

The proof of Conjecture~\ref{conj-b} for type \(A_{n-1}\) is a
rather long calculation.  We were
motivated and guided by the calculations in \cite{KM}, where the case
\(n=3\) is treated.  The key is to establish a length-preserving
decomposition of the form
\[
W=X_0 \times \bbT_{1} \times X_{1} \times \cdots \times
 \bbT_{n-1}\times X_{n-1},
\]
where \(\bbT_i=\{w_{i}^N:N \geq 0\}\) and \(X_i\) are certain finite
sets.

Define \(\scrS_n\) as the closed fundamental alcove and set \(\scrF_0
= \scrA\). A collection of subsets of \(\scrA\) is deemed almost
disjoint if any two distinct members intersect in measure zero. The above decomposition requires \(\scrF_1, \dots, \scrF_{n-1}\) and \(\scrS_1, \dots, \scrS_{n-1}\) satisfying:
\begin{description}
\item[Claim 1.] For every \(0 \leq i \leq n-1\), \(\scrF_{i} =
  \bigcup_{w \in X_{i}} w\scrS_{i+1}\) and \(\{w\scrS_{i+1}\}_{w \in
    X_i}\) is almost disjoint.
\item[Claim 2.] For every \(1 \leq i \leq n-1\), \(\scrS_{i} =
  \bigcup_{w \in \bbT_{i}} w\scrF_{i}\) and \(\{w\scrF_i\}_{w\in
    \bbT_i}\) is almost disjoint.
\end{description}
Thus 
\[
\scrF_0\supset\scrS_1\supset \scrF_1\supset\scrS_2\supset\cdots\supset\scrS_{n-1}\supset
\scrF_{n-1}\supset \scrS_n.
\]

The calculations in \cite{KM} suggest the following choice of region:
\begin{align*}
    \scrS_i &:= \{ x_0 \leq  \cdots  \leq  x_{i-1} \leq \min\{x_{i}, x_{i+1}, \ldots, x_n\} \}.
\end{align*}
We define \(\scrS_i\) by this formula for
\(i=1,\ldots,n\).  Observe that \(\scrS_n\) is indeed the
closed fundamental alcove.

One verifies that \(\scrS_{i+1} \subset \scrS_i\),
\(w_{i}(\scrS_i) \subseteq \scrS_i\), and \(\bigcap_{N \geq 0}
w_{i}^N(\scrS_i)\) is empty.   As a result, the set
\[
    \scrF_i := \overline{\scrS_i \setminus w_{i}(\scrS_i)}
    \qquad (1 \leq i \leq n-1)
\]
is a fundamental domain for the monoid \(\bbT_i\) acting on \(\scrS_i\).  In other words, with these definitions of
\(\scrS_i\) and \(\scrF_i\), Claim 2 is already established for
\(i=1,\ldots,n-1\).

\subsection{The sets \(X_i\)}

It remains to construct the finite sets \(X_i\) giving the
complementary pieces in Claim 1.

The proof of Claim 1 is divided into three cases.

\begin{theorem} \label{Xi}
  For \(0\leq i\leq n-1\), there is unique subset \(X_i\) of \(W\)
  such that
  \begin{itemize}
    \item[{\rm (a)}] \(\scrF_{i} =
  \bigcup_{w \in X_{i}} w\scrS_{i+1}\) and \(\{w\scrS_{i+1}\}_{w \in
    X_i}\) is almost disjoint;
  \item[{\rm(b)}] \(p_{X_i}=(1-u^n)/(1-u)\); in other words,
    \(|X_i|=n\) and \(X_i\) has a unique element of length \(j\) for \(j=0,\ldots,n-1\).
  \end{itemize}
  In fact, \(X_i=\{\sigma_0,\ldots,\sigma_{n-1}\}\) is a subset of
  \(\Waff_{v_i}\) with
\[ 
\sigma_{j} = \begin{cases}
(x_i,\ldots, x_{i+j})^{-1} &\text{if } j \leq n-i, \\
(x_i,\ldots, x_n)^{-1}(x_{n}, \ldots,x_{i+j})  &\text{if } j \ge  n-i.
\end{cases}
\]
\end{theorem}

Observe that  
\[
  \sigma_j=\begin{cases}
    s_{i+j}\cdots s_{i+1} &\text{ if } j\leq n-i,\\
    (s_n\cdots s_{i+1})(s_{n+1}\cdots s_{i+j})&\text{ if } j\geq n-i
  \end{cases}
\]
is a product of \(j\) distinct simple reflections.  Therefore, \(\ell(\sigma
_j)=j\) and (b) is true.  The rest of this section is devoted to the
proof of (a) when \(X_i\) is specified by the theorem.  The uniqueness
of \(X_i\) will be proved in the next section.

We remark that in the description of \(\sigma_j\) above, the
convention given at the end of \ref{AnSetting} is in action.  This is
very useful for computing \(\sigma _j(\scrR)\) when \(\scrR\) is a
region defined by inequalities given in terms of \(x_i,\ldots,x_{i+n-1}\).
Then \(\sigma _j(\scrR)\) is the region defined by the same
inequalities transformed by \(\sigma _j\).

\subsubsection{The case \(i=0\)}
Given that 
\[
\scrS_{1}=\{(x_0,\ldots,x_{n-1}): x_0 \text{ is minimal among } x_0,\ldots,x_{n-1}\}
\]
and 
\( 
\sigma_j\scrS_{1}=\{(x_0,\ldots,x_{n-1}): x_j \text{ is minimal among } x_0,\ldots,x_{n-1}\}
\)
it follows that 
\( 
\bigcup_{j=0}^{n-1}\sigma_j \scrS_{1} = \scrA = \scrF_0
\) 
forms an almost disjoint union.

\subsubsection{The case \(i=n-1\)}

Writing in terms of \(x_{n-1},\ldots,x_{2n-2}\), which is convenient
for working with \(\Waff_{v_{n-1}}\),
we have \(\scrS_n=\{x_{n-1}\leq x_n\leq \cdots\leq x_{2n-2}\leq
x_{n-1}+1\}\).  Also, we have
\(\sigma_j=(x_{n-1},x_n,\ldots,x_{n+j-1})\).  Therefore,
\begin{align*}
\sigma_j\scrS_n =&\{x_n\leq \cdots\leq
x_{2n-2}\leq x_n+1\} \\
&\cap\left\{\begin{array}{ll}
                       x_{2n-2}-1\leq x_{n-1}\leq x_n &\mbox{ if }j=0\\
                       x_{n+j-1}\leq x_{n-1}\leq x_{n+j}&\mbox{ if }
                       j=1,\ldots,n-2\\
                       x_{2n-2}\leq x_{n-1}\leq x_n+1 &\mbox{ if }j=n-1\end{array}\right\}
\end{align*}
Thus \(\bigcup_{j=0}^{n-1}\sigma_j \scrS_n=\{x_n\leq \cdots\leq
x_{2n-2}\leq x_n+1\}\cap\{x_{2n-2}-1\leq x_{n-1}\leq x_n+1\}\).

On the other hand, \[\scrS_{n-1}=\{x_n\leq \cdots\leq
x_{2n-2}\leq x_n+1\}\cap\{x_{2n-2}\leq x_{n-1}+1\}\] 
and 
\[w_{n-1}\scrS_{n-1}=\{x_n\leq \cdots\leq
x_{2n-2}\leq x_n+1\}\cap\{x_n+1\leq x_{n-1}\}.\]  So
\begin{align*}
\scrF_{n-1} &= \overline{\scrS_{n-1}\setminus w_{n-1}\scrS_{n-1}}\\
&=\{x_n\leq \cdots\leq x_{2n-2}\leq x_n+1\}\cap\{x_{2n-2}-1\leq x_{n-1}\leq x_n+1\}.
\end{align*}
 This
proves the desired result that \(\scrF_{n-1} = \bigcup_{j=0}^{n-1}\sigma_j \scrS_{n}\) which is an almost disjoint union.

\subsubsection{The case \(1 \leq i \leq n-2\)}
Fix an \(i\) in the range \(1 \leq i \leq n-2\).  We have
\(\scrS_i=\{x_n\leq \cdots\leq
x_{n+i-1}\leq\min\{x_i,\ldots,x_n\}+1\}\) and
\[
  \scrS_{i+1}=\scrS_i\cap\{x_{n+i-1}-1\leq \min\{x_i,\ldots,x_n\}=x_i\}.
\]
It follows that for \(0\leq j \leq n-i-1\), we have \(\sigma
_j\scrS_i=\scrS_i\) and 
\[
\sigma _j\scrS_{i+1}=\scrS_i\cap\{x_{n+i-1}-1\leq \min \{x_i,\ldots,x_n\}=x_{i+j}\}.
\]
It follows that the collection \(\{\sigma_j\scrS_{i+1}\}_{0\leq j\leq n-i-1}\) is
almost disjoint and
\begin{align*}
\scrF_i'&:=\bigcup_{j=0}^{n-i-1} \sigma_j \scrS_{i+1} = \scrS_i \cap \left\{
  x_{n+i-1} -1\leq \min\{ x_{i}, \ldots,x_n \} \leq x_n \right\}\\
&=\scrS_i\cap\{\min\{ x_{i}, \ldots,x_{n-1} \} \leq x_n \}.
\end{align*}
The region \(\scrF_i'\)  is almost disjoint from
\[
w_i(\scrS_i) = \scrS_i \cap \left\{ x_n \leq \min\{ x_{i}, \ldots, x_{n-2}, x_{n-1}-1 \} \right\},
\]
and we have
\[@
\overline{\scrS_i\setminus \scrF_i'}&=\scrS_i \cap \left\{ x_n \leq \min\{ x_{i},
  \ldots, x_{n-1} \} \right\},\\
\scrF_i'':=\overline{(\scrS_i\setminus \scrF_i')\setminus
  w_i(\scrS_i)}&=\scrS_i\cap\{x_{n-1}-1\leq x_n\leq \min\{x_i,\ldots,x_{n-1}\}\}.
\]
Thus, \(\scrF_i=\scrF_i'\cup\scrF_i''\) and \(\scrF_i'\) and
\(\scrF_i''\) are almost disjoint.
It remains to show that the collection \(\{\sigma_j\scrS_{i+1}\}_{n-i\leq j\leq
  n-1}\) is almost disjoint and their union is \(\scrF_i''\).

Now assume \(n-i \leq j\leq n-1\).  Apply
\(\sigma_j=(x_n,\ldots,x_{i+j},x_{n-1},x_{n-2},\ldots,x_i)\) to
\(\scrS_{i+1}=\{x_n\leq \cdots\leq x_{n+i-1}\leq
x_i+1\}\cap\{x_i\leq \min\{x_{i+1},\ldots,x_n\}\}\),we have that \(\sigma_j\scrS_{i+1}\) is the intersection of \(\{x_n\leq\min\{x_i,\ldots,x_{n-2},x_{n+1}\}\}\) and 
\[\left\{\begin{array}{ll}
x_{n+1}\leq \cdots\leq x_{i+j}\leq x_{n-1}\leq
x_{i+j+1}\leq \cdots\leq x_{n+i-1}\leq
                                 x_n+1&\mbox{if } j\leq n-2\\
x_{n+1}\leq \cdots\leq x_{n+i-1}\leq x_{n-1}\leq x_n+1&\mbox{if } j=n-1\end{array}\right\}\]

Observe that this is contained in \(\scrS_i\) (for example, to see
that every point on \(\sigma _j\scrS_{i+1}\) satisfy \(x_{n+i-1}\leq
x_{n-1}+1\), use \(x_{n+i-1}\leq x_n+1\leq x_{n+1}+1\) and
\(x_{n+1}\leq x_{n-1}\)).  We then write
\begin{align*}
  \sigma _j\scrS_{i+1}=&\scrS_i\cap\left\{
    \begin{array}{ll}
      x_{i+j} \leq x_{n-1}\leq x_{i+j+1} &\mbox{ if } j\leq n-2\\
      x_{n+i-1}\leq x_{n-1}\leq x_n+1 &\mbox{ if } j=n-1
    \end{array}\right\}\\
    &\cap\{x_n\leq \min\{x_i,\ldots,x_{n-2}\}\}.
\end{align*}
This makes it clear that the collection \(\{\sigma_j\scrS_{i+1}\}_{n-i\leq j\leq
  n-1}\) is almost disjoint and their union is
\[
\scrS_i\cap\{x_n\leq x_{n-1}\leq x_n+1, x_n\leq \min\{x_i,\ldots,x_{n-2}\}\},
\]
which coincides with \(\scrF_i''\).

\subsection{Length-preserving decomposition}

The almost-disjoint paving now translates into a product decomposition
of the affine Weyl group; the remaining point is length preservation.

Assume that \(\{X_i\}_{1\leq i\leq n-1}\) are sets satisfying
condition (a) of Theorem~\ref{Xi}.  Then
\(W=X_0\times \bbT_1\times
X_1\times \cdots \times\bbT_{n-1}\times X_{n-1}\) is a decomposition.  It remains to show that
if condition (b) of Theorem~\ref{Xi} is also satisfied, then the
decomposition is length-preserving.

The Poincar\'e series of the affine Weyl group of type \(\tilde{A}_{n-1}\) can be deduced via Bott's formula in Section 8.9 of \cite{Hum}:
\(
p_W = (1 - u^n)/(1 - u)^n.\)
Together with \(p_{\bbT_i}=1/(1-u^n)\) for all \(i\) and condition
(b), we deduce
\[
 p_W =  p_{X_0} p_{\bbT_{1}} p_{X_{1}} \cdots p_{\bbT_{n-1}} p_{X_{n-1}}.
\]
Invoking Proposition~\ref{lengthprevseringdecompsition}, we obtain:
\begin{equation}
   P_{W} =  P_{X_0} P_{\bbT_{1}} P_{X_{1}} \cdots P_{\bbT_{n-1}} P_{X_{n-1}}  \label{eq:An1}
\end{equation}
and  the length-preserving decomposition:
\[
    W = X_0 \times \bbT_{1} \times X_{1} \times \cdots \times \bbT_{n-1} \times X_{n-1}.\label{decomp:An1}
\]
We next prove the uniqueness part of Theorem~\ref{Xi}.  Suppose
that \(\{X_i\}_{1\leq i\leq n-1}\) and \(\{X_i'\}_{1\leq i\leq n-1}\)
both fulfill the required conditions of Theorem~\ref{Xi}.  We
prove \(X_i=X_i'\) by induction on \(i\).  The case \(i=n-1\)
follows directly from the fact that \(\scrS_n\) is a fundamental
domain for the action of \(W\) on \(\scrA\).  Assume that we have
\(X_i=X_i'\) for \(i=i_0+1,\ldots,n-1\).  Put \(L=\bbT_{i_0+1}\times
X_{i_0+1}\times\cdots\times \bbT_{n-1}\times X_{n-1}\).  Then we have
\[
X_{i_0}\times L=X_{i_0}'\times L=\{w\in W: w\scrS_n\subset
\scrF_{i_0}\}
\]
and the decompositions \(X_{i_0}\times L\) and \(X_{i_0}'\times L\)
are length-preserving.  By Proposition~\ref{uniq-dec}, we have \(X_{i_0}=X_{i_0}'\).

\subsection{The Proof of Conjecture \ref{conj-b} for Type $\tilde{A}_{n-1}$ }\label{conjAn}

With the monoid factorization established, we combine it with the
parabolic reduction from Proposition~\ref{irrCoxeter} to obtain the
required identity in the abelian quotient.

We prove Conjecture \ref{conj-b}, which claims that 
\[
\prod_{I \subseteq S} P_{W_I}^{(-1)^{|S \setminus I|}} = \prod_{i=1}^{n-1} P_{\bbT_i} \quad \text{in} \, \left(H[[u]]^\times\right)^{\text{ab}}.
\]
By Proposition \ref{irrCoxeter}, we have
\begin{equation} \label{eq:An2}
\prod_{I \subseteq S} P_{W_I}\!^{(-1)^{|S \setminus I|}} = \prod_{J:\text{irreducible, dense}} P_{W_J}\!^{(-1)^{|S \setminus J|}} \quad \text{in} \, \left(H[[u]]^\times\right)^{\text{ab}}.
\end{equation}
The proper irreducible dense subsets \(J\) of \(S\) are:
\begin{itemize}
\item For size $n-1$: $J_i = \{s_{i+1}, \cdots, s_{n+i-1}\}$ for $i=0$ to $n-1$.
\item For size $n-2$: $J'_i = \{s_{i+2}, \cdots, s_{n+i-1}\}$ for $i=0$ to $n-1$.
\end{itemize}
Combining Equations (\ref{eq:An1}) and (\ref{eq:An2}), we have, in \(\left(H[[u]]^\times\right)^{\text{ab}}\),
\[@
\prod_{I \subseteq S} P_{W_I}\!^{(-1)^{|S \setminus I|}} 
&=P_{\Waff} \prod_{i=0}^{n-1} P_{W_{J_i}}\!^{(-1)} P_{W_{J_i'}}\\
& = \left(\prod_{i=1}^{n-1} P_{\bbT_i}\right)\left( \prod_{i=0}^{n-1} P_{X_{i}} P_{W_{J_i}}\!^{(-1)} P_{W_{J_i'}}  \right).
\]
It remains to show that for $0 \leq i \leq n-1$,
\[ P_{X_{i}}  = P_{W_{J_i}} P_{W_{J_i'}}\!^{(-1)} \quad \text{in} \,
  \left(H[[u]]^\times\right)^{\text{ab}},
\]
which follows from Theorem~\ref{idenAn}.





\section{Result for Type \(\tilde{C}_{n}\)}

The type \(\tilde C_n\) case follows the same strategy as type
\(\tilde A_{n-1}\), but only two auxiliary finite factors are needed.

\subsection{Setup}
For type \(\tilde{C}_{n}\), take \(\scrA=\bbR^{n}\). For \(\vec{\xi} = (\xi_1, \ldots, \xi_{n}) \in \scrA\), set \( x_i(\vec{\xi}) = \xi_i \). The affine functions
\[
\{ a_0 := 1 - 2 x_1, a_1 := x_1 - x_2, \ldots, a_{n-1} := x_{n-1} - x_{n}, a_{n} := 2 x_{n} \}
\]
forms a system of simple affine roots.

The corresponding closed fundamental alcove  is defined by
\[
 \left\{ \frac{1}{2} \geq x_1 \geq \cdots \geq x_{n} \geq 0 \right\}.
\]
The fundamental coweights are
\[
\omega_i^\vee = (\underbrace{1, 1, \ldots, 1}_{i\text{-times}}, 0,
\ldots, 0) \quad (1\leq i \leq n-1),\quad\mbox{and}\quad \omega_{n}^\vee = \frac{1}{2} (1, \ldots, 1).
\]
The sum of positive roots is
\[
2 \rho = 2 n x_1 + 2 (n-1) x_2 + \cdots + 4 x_{n-1} + 2 x_{n}.
\]
The roots are \( \pm x_i \pm x_j \) for all \(1 \leq i < j \leq n\) and \( \pm 2 x_i \) for \(1 \leq i \leq n\).

We endow \(\bbR^{n}\) with the standard inner product.  Let \( s_i
\) denote the orthogonal reflection through the hyperplane \(\{a_i = 0
\}\) and put \(S=\{s_0,\ldots,s_{n}\}\), \(W=\langle s_0,\ldots,s_{n}\rangle\).

\subsubsection{Geodesic tube with straight generator}
Consider
\[
w_{n}(\vec{\xi}) = s_0 s_1 \cdots s_{n}(\vec{\xi}) = (\xi_{n} + 1, \xi_1, \ldots, \xi_{n-1})
\]
and for \(1 \leq i \leq n-1\),
\[
w_i(\vec{\xi}) = \left( w_{n} s_{n-1} \cdots s_i \right)(\vec{\xi}) = (\xi_i + 1, \xi_1, \xi_2, \ldots, \xi_{i-1}, \xi_{i+1}, \xi_{i+2}, \ldots, \xi_{n}).
\]
Then $w_i$ has a word expression of length \(2n-i+1\).
By direct computation, we obtain
\[
(w_i)^i(\vec{x}) = \vec{x} + \begin{cases}
\omega_i^\vee, & \text{if } i \leq n-1; \\
2\omega_i^\vee, & \text{if } i = n.
\end{cases}
\]
This implies that \( w_i \) is hyperbolic with \( v_{w_i} = \frac{1}{i} \omega_i^\vee \) for \( i = 1 \) to \( n-1 \) and \( v_{w_{n}} = \frac{2}{n} \omega_{n}^\vee \).

The set
\[
\scrT_i:=\{x_i+1 > x_1 > \cdots > x_{i}, \frac{1}{2}>x_{i+1}> \cdots > x_{n} > 0 \}
\]
is a geodesic tube of direction \(\omega_i^\vee\), and direct verification shows that \(\scrT_i\) is \(w_i\)-invariant.  Hence \(w_i\in \Stab^h(\scrT_i)\). Here \(a_i^\vee\) denotes the coroot corresponding to the finite simple root underlying the affine simple root \(a_i\). From the
defining formula, we have \(\ell(w_i)\leq 2n-i+1\).  By
Theorem~\ref{straightlength}, \(\ell(w_i)\geq
\frac{(\omega_i^\vee,a_i^\vee)}{(\omega_i^\vee,\omega_i^\vee)}\langle
2\rho,\omega_i^\vee\rangle =2n-i+1\).  Thus we conclude that \(w_i\) is the
canonical generator of \(\scrT_i\), and it is straight.

\subsection{The Regions $\scrS_i$}

We now choose nested regions on which the monoids generated by the
\(w_i\) act with explicit fundamental domains.

Let \(\scrF_0 =\scrA\) and
\[
\scrF_{n}=\{ 1  \geq x_1\geq \cdots \geq x_{n} \geq 0\}.
\]
Let $\scrS_{n+1}$ be the closed fundamental alcove and for $1 \leq i \leq n$, let
\[
\scrS_{i}=\{x_{i}+1 \geq x_1\geq \cdots \geq x_{n} \geq 0\}.
\]
Then for \(1\leq i\leq n\), we have
\begin{align*}
w_i(\scrS_i)&=\{x_1 \geq x_2 \geq \cdots \geq x_i \geq x_1-1 \geq x_{i+1} \geq \cdots \geq x_{n} \geq 0\}\\
&=\{x_1 \geq  \cdots \geq x_{n} \geq 0, x_i+1 \geq x_1 \geq x_{i+1}+1\} \subset \scrS_i.
\end{align*}
The formula above shows that $\bigcap_{N\geq 0}
w_i^N(\scrS_i)$ is empty, and
\(\scrS_{i+1}=\overline{{\scrS_i\setminus w_i(\scrS_i)}}\) can be
regarded as a fundamental domain for the monoid $\bbT_i=\{w_i^N: N\geq
0\}$ acting on \(\scrS_i\), for \(1\leq i\leq n-1\).  Similarly, \(\scrF_{n}\) is a fundamental
domain for the monoid \(\bbT_n=\{w_n^N:N\geq 0\}\) acting on \(\scrS_n\).

Note that we have 
\[ \scrA = \scrF_0 \supset \scrS_1 \supset \cdots \supset \scrS_{n} \supset \scrF_{n} \supset \scrS_{n+1}.\]
We construct two subsets \(X_0\) and \(X_n\) of \(\Waff\) such that \(\scrF_i = \bigcup_{ g \in X_i} g \scrS_{i+1}\) is almost disjoint for \(i=0,n\), and such that
\[
 X_0 \times  \bbT_{1} \times \cdots \times \bbT_{n} \times X_{n}
\] 
is a length-preserving decomposition of \(\Waff\), as in type \(\tilde{A}_{n-1}\).

\subsection{The sets \(X_0\) and \(X_{n}\)}

The only remaining pieces in the type \(\tilde C_n\) paving are the
initial and terminal finite factors.

Let \(v_i\) be the vertex of the closed fundamental alcove fixed by all elements in \(S\) except \(s_i\). Let \(\Waff_{v_i}\) be the stabilizer of \(v_i\), the parabolic subgroup generated by \(S \setminus \{s_i\}\), and set \(\Waff_{v_0,v_{n}} = \Waff_{v_0} \cap \Waff_{v_{n}}\).

We set \(X_0= \Waff_{v_0}\) and
\begin{align*} 
X_{n} &= \Waff_{v_0,v_{n}} \backslash \Waff_{v_{n}} \\
&= \{w \in  \Waff_{v_{n}} : \ell(gw) = \ell(w) + \ell(g), \forall g \in \Waff_{v_0,v_{n}}\} \\
&= \{w \in  \Waff_{v_{n}} : \ell(s_iw) = \ell(w) + 1, \text{ for } i = 1, \ldots, n-1 \}
\end{align*}

\begin{proposition} \label{step1}
The set \(\scrF_i = \bigcup_{ g \in X_i} g \scrS_{i+1}\) is an almost disjoint union  for $i=0$ and $n$. 
\end{proposition}

\begin{proof}
  Let $i=0$.  Then \(X_0=\Waff_{v_0}\) is the Weyl group of type
  \(C_{n}\) and $S_{1} = \{ x_1\geq \cdots \geq x_{n} \geq
  0\}$ is a fundamental domain for the action of \(X_0\) on \(\scrA\)
  (see \cite[Planche III]{Bou}).  This proves the result in this case.

Let $i=n$. 
Here
\[
\scrS_{n+1} = \left\{ \frac{1}{2} \geq x_1 \geq \cdots \geq x_{n} \geq 0 \right\}.
\]
Since \(\Waff_{v_0,v_{n}}\) includes all permutations on the coordinates, we have
\[
\bigcup_{g \in \Waff_{v_0,v_{n}}} g\scrS_{n+1}  = \left\{\frac{1}{2} \geq x_1, \cdots, x_{n} \geq 0\right\}.
\]
Combined with the action of \(s_0\) mapping \(x_1\) to \(1 - x_1\), we get
\[
\bigcup_{g \in \Waff_{v_0}} g\scrS_{n+1}  = \{1 \geq x_1, \cdots, x_{n} \geq 0\}.
\]

For any \(\vec{\xi} \in \scrS_{n+1}\) and \(w \in \Waff\), we have
\[
\ell(s_iw) = \ell(w) + 1 \iff a_i(w(\vec{\xi})) = (w\xi)_i - (w\xi)_{i+1}\geq 0.
\]
Therefore,
\begin{align*}
\bigcup_{g \in X_{n}} g\scrS_{n+1} &= \{1 \geq x_1, \cdots, x_{n} \geq 0\} \cap \left(\bigcap_{i=1}^{n-1} \{x_i - x_{i+1} \geq 0\}\right) \\
&= \{1 \geq x_1 \geq x_2 \geq \cdots \geq x_{n} \geq 0\} = \scrF_{n}.
\end{align*}
As \(\scrS_{n+1}\) is the closed fundamental alcove, the above union is almost disjoint.
\end{proof}

\subsection{The Proof of Conjecture \ref{conj-b} for Type $\tilde{C}_{n}$ }

The preceding paving gives the required length-preserving
decomposition, so the final step is to compare its Poincar\'e series
with the alternating parabolic product.

The preceding construction gives a decomposition
\begin{equation}\label{Cndecomp}
\Waff = X_0 \times \bbT_1 \times \cdots \times \bbT_{n} \times X_{n}.
\end{equation}
Here \(X_0=\Waff_{v_0}\) and  \(X_{n}=\Waff_{v_0,v_{n}}\backslash \Waff_{v_{n}}\),
whose Poincar\'e series are computed by Chevalley's formula in Section 3.15 of \cite{Hum}.  Bott's formula in Section 8.9 of \cite{Hum} gives the Poincar\'e series of the affine Weyl group of type $\tilde{C}_{n}$, and hence
\[
p_W =  p_{X_0} p_{\bbT_{1}} \cdots p_{\bbT_{n}} p_{X_{n}}.
\]
Invoking Proposition~\ref{lengthprevseringdecompsition}, we have:
\[
   P_{W} =  P_{X_0} P_{\bbT_{1}} \cdots P_{\bbT_{n}} P_{X_{n}} 
  = P_{\Waff_{v_0}} P_{\bbT_{1}} \cdots P_{\bbT_{n}} (P_{\Waff_{v_0,v_{n}}}\!)^{-1} P_{\Waff_{v_{n}}}
\]
and  the decomposition \ref{Cndecomp} is length-preserving. 

Applying Proposition \ref{irrCoxeter}, note that the irreducible dense subsets are $S$, $S\setminus\{s_0\}$,  $S\setminus\{s_{n}\}$, and  $S\setminus\{s_0,s_{n}\}$.
Thus, in $\left(H[[u]]^\times\right)^{\text{ab}}$,
\begin{align*}
\prod_{I \subset S} P_{W_I}\!^{(-1)^{|S \setminus I|}} &= \prod_{J:\rm{irreducible, dense}} P_{W_J}\!^{(-1)^{|S \setminus J|}} \\
&= P_{\Waff}(P_{\Waff_{v_0}})^{-1}(P_{\Waff_{v_{n}}})^{-1}P_{\Waff_{v_0,v_{n}}} = \prod_{i=1}^{n} P_{\bbT_i}.
\end{align*}

This completes the proof of Conjecture \ref{conj-b}.

\begin{remark} (i) The sets  \(X_0\) and \(X_{n}\) are uniquely determined in a
  manner similar to Theorem~\ref{Xi}.  (ii)  It is
possible to prove (\ref{Cndecomp}) by showing
\(W=X_{n}^{-1}\times \bbT_{n}^{-1}\times \cdots \times
\bbT_1^{-1}\times X_0^{-1}\).  This involves paving up \(\scrA\) in a
different way, and the role of \(\scrS_i\) will be played by \((-\infty,1/2]^{n+1-i}\times[-1/2,1/2]^{i-1}\).  We leave out the details.
\end{remark}

\end{document}